\documentclass[3p,number,preprint,times]{elsarticle}
\usepackage{graphicx}
\usepackage{subcaption}
\usepackage{siunitx} 
\usepackage{multirow}
\usepackage{xcolor}
\usepackage{longtable}
\usepackage{hhline}
\usepackage[linesnumbered,ruled]{algorithm2e}
\usepackage{amsmath}
\usepackage{lineno}
\usepackage{parskip}
\usepackage{makecell}

\graphicspath{ {./images/} }

\journal{EURO Journal on Transportation and Logistics}

\biboptions{authoryear}

\begin{document}
	
\begin{abstract}
		
		Recent developments in modular transport vehicles allow deploying multi-purpose vehicles which can alternately transport different kinds of flows. In this study, we propose a novel variant of the pickup and delivery problem, the multi-purpose pickup and delivery problem, where multi-purpose vehicles are assigned to serve a multi-commodity flow. We solve a series of use case scenarios using an exact optimization algorithm and an adaptive large neighborhood search algorithm. We compare the performance of a multi-purpose vehicle fleet to a mixed single-use vehicle fleet. Our findings suggest that total costs can be reduced by an average of 13\% when multi-purpose vehicles are deployed, while at the same time reducing the total vehicle trip duration and total distance travelled by an average of 33\% and 16\%, respectively. The size of the fleet can be reduced by an average of 35\%. The results can be used by practitioners and policymakers to decide on whether the combination of passenger and freight demand flows with multi-purpose vehicles in a given system will yield benefits compared to existing fleet configurations.
		
		\fbox{\parbox{0.9\textwidth}{Highlights:
				\begin{itemize}
					\item formulating and solving a novel multi-purpose pickup and delivery problem
					\item analysis of urban scenarios with varied spatial and temporal demand
					\item cost savings by an average of 13\% 
					\item reduction of fleet size by an average of 35\% 
					\item reduction of total vehicle trip duration by an average of 33\%
					\item reduction of total distance travelled by an average of 16\%
				\end{itemize}
			}
		}
\end{abstract}
	
\begin{keyword}
		Public transportation \sep
		Freight transportation \sep
		Multi-purpose vehicles \sep
		Heuristic optimization 
\end{keyword}

\begin{frontmatter}
		\title{Multi-purpose Pickup and Delivery Problem for Combined Passenger and Freight Transport 
		}

		\author{Jonas Hatzenb{\"u}hler\fnref{myfootnote}\corref{mycorrespondingauthor}}
		\ead{jonashat@kth.se}
		\author{Erik Jenelius\fnref{myfootnote}}
		\author{Gy\H{o}z\H{o} Gid\'ofalvi \fnref{secondfootnote}}
		\author{Oded Cats \fnref{myfootnote,thirdfootnote}}
		
		\cortext[mycorrespondingauthor]{Corresponding author:}

		\address{KTH Royal Institute of Technology, 100 44 Stockholm, Sweden}
		
		\fntext[myfootnote]{Division of Transportation Planning, KTH Stockholm, Sweden}
		\fntext[secondfootnote]{Department of Urban Planning and Environment, KTH Stockholm, Sweden}
		\fntext[thirdfootnote]{Department of Transport and Planning, Delft University of Technology, Netherlands}
		
		\date{Submitted: June, 2022}
		
\end{frontmatter}

\section{Introduction}\label{sec:introduction}

The ongoing trends of urbanization, increasing e-shopping, and digitalization of supply chain management challenge the efficiency and sustainability of the existing transportation systems. \citet{savelsbergh50thAnniversaryInvited2016} and \citet{losValueInformationSharing2020} stress the importance of integrated and/or collaborative transportation solutions to achieve efficiency and sustainability. The authors specifically address the concept of utilizing vehicle capacities for multiple purposes.

Public transportation and urban freight delivery operations are conventionally designed as two independent systems. Public transportation systems are designed to satisfy the passenger demand peaks \citep{cederBusNetworkDesign1986}, while logistic operations are designed to meet delivery times \citep{ghilasPickupDeliveryProblem2016}. In off-peak periods public transportation may experience unused vehicle capacity and empty-vehicle kilometers, while the urban logistics during the same periods are dominated by daily goods and parcel deliveries between stores, warehouses, and companies. For traffic safety reasons and maneuverability in urban environments, these vehicles are typically small or medium-sized trucks. However, both transportation systems strive to efficiently satisfy the transportation needs of requests (i.e., passenger or freight items) at a specific time to a specific destination. 

The concept of combining multiple demand flows in one transport system is known as integrated transportation. In the past, several projects have investigated and demonstrated the successful integration of passenger and freight transportation in urban environments. In Amsterdam, Netherlands \citep{marinovUrbanFreightMovement2013} a pilot project investigated the delivery of consumer goods in specially designed trams that operate in between the passenger trams. The trams transported goods from a suburban depot to a shopping mall in the city center. The integration of flows was achieved by sharing the same track network and infrastructure. The service was operating from 2007-2009 but had to be terminated due to missing political and financial support since the operations could not be realized without subsidies. A vehicle production plant in Dresden, Germany \citep{article} was supplied by a tram. The tram operates during off-peak hours and uses the same tracks as the passenger tram. The tram has only two stops, at the production plant and at the logistic center. The operations in Dresden were terminated in the beginning of 2021 after 19 years of operations. The reason for this was a changed logistic concept for the new vehicle produced in the factory, which fully utilized delivery trucks.

The consolidation of multiple demand flows in one transport system can be realized in different ways. \citet{mouradSurveyModelsAlgorithms2019} present a survey of models and algorithms for optimizing the shared transportation of passenger and freight items. 

One concept for integrating multiple item types is sequential integration. In this concept, different types of items are transported by the same vehicle but at different times. \citet{steadieseifiMultimodalFreightTransportation2014}, \citet{cochraneMovingFreightPublic2017}, \citet{ozturkOptimizationModelFreight2018}, \citet{behiriUrbanFreightTransport2018a} and \citet{liUrbanRailService2021} propose different modes of operation and give examples of successful integration. This concept is typically applied to rail-bound transportation modes like the tram, metro, or train systems or to systems combining fixed-line services with on-demand services. In \cite{mouradIntegratingAutonomousDelivery2020} the authors propose a variant of the PDP with simultaneous integration. Small pickup and delivery robots are integrated into scheduled line services. The robots can then be used to directly serve certain requests in a local area or travel in the scheduled line to reach other requests. The authors show that this integrated approach can lead to 18.2\% cost savings compared to a pure freight system. Compared to simultaneous integration concepts, sequential integration requires to physically adjust the vehicle for each item type. However, the combined optimization of multiple item types may still lead to an overall improved system as compared to the independent operation of single-purpose fleets. Additionally, for sequential integration of multiple item types, passenger requests are not mixed with e.g. freight items, which guarantees that passenger requests are not affected by the delivery and pickup of freight items. For simultaneous integration concepts, freight requests would be mixed with passenger requests which would lead to additional stops for the vehicles, which in turn reduce the level of service for the passengers in this vehicle. Another advantage of sequential integration is the higher capacity for each demand type when exclusively served, which brings greater flexibility for which freight items can be delivered.

The efficient operation of multiple vehicles is first already envisioned and described in \cite{dantzigTruckDispatchingProblem1959} as the truck dispatching problem. Since then numerous variations, extensions, and solution algorithms have been developed (see \cite{tothVehicleRoutingProblems2014, kocThirtyYearsHeterogeneous2016, dundarReviewSustainableUrban2021} for recent comprehensive literature reviews). Two of the most studied variations are the Pickup and Delivery Problem (PDP), where a request has a pickup position and a delivery position, and the more general Dial-a-ride Problem (DARP), in which dynamic requests are served (see \cite{PDP1988}, \cite{savelsberghGeneralPickupDelivery1995} for the PDP, and \cite{cordeauDialarideProblemModels2007}, \cite{berbegliaDynamicPickupDelivery2010} for the DARP).

\begin{figure*}[ht]
	\centering
	\begin{subfigure}[b]{0.49\textwidth}
		\centering
		\includegraphics[width=\textwidth]{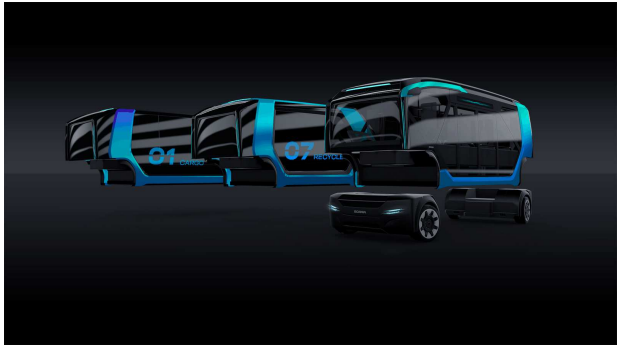}
		\caption[]%
		{Scania NXT vehicle concept \citep{scaniaNXTConceptVehicle2020}}
		\label{fig:nxt}
	\end{subfigure}
	\hfill    
	\begin{subfigure}[b]{0.49\textwidth}
		\centering
		\includegraphics[width=\textwidth]{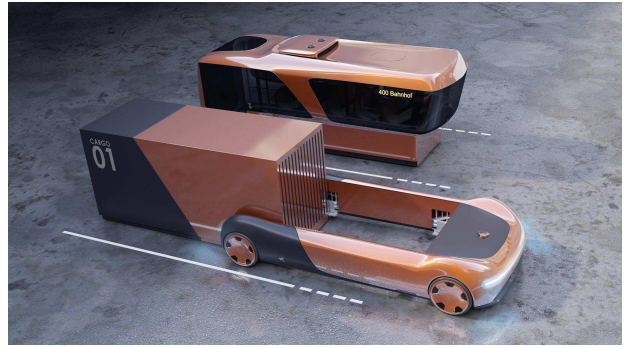}
		\caption[]%
		{U-Shift vehicle concept \citep{UShiftIIDemonstrator}}
		\label{fig:ushift}
	\end{subfigure}
	\caption{Illustration of the multi-purpose vehicle concept}
	\label{fig:illustration}
\end{figure*}

In this study, we revise the pickup and delivery problem formulation to facilitate the operation of multi-purpose vehicles, i.e. the Multi-Purpose Pickup and Delivery Problem (MP-PDP). We analyze in a series of experiments the potential benefits this vehicle concept has over conventional transportation systems. The proposed problem is a variation of the heterogeneous routing problem (see \cite{kocThirtyYearsHeterogeneous2016}) and the truck-and-trailer routing problem (see \cite{derigsTruckTrailerRouting2013}). The main variation in comparison to these problems is the change of vehicle type during a vehicle trip. Additionally, the module can be changed between several vehicle platforms at several places, whereas at the truck-and-trailer routing problem, the trailer is dropped and connected to the truck at the same dedicated place. In the following, we study the operation of sequential deliveries of freight and passengers between origins (pickup locations) and destinations (drop off locations). Hence, our work also extends the sequential integration concepts for road network-based operations by adding the concept of modularity and vehicle module handling to the problem formulation.

The operation can be explained as follows: A vehicle operates on a route that starts and ends at a depot. A vehicle can either deliver goods from the depot to a set of customers or transport passengers from their origins to their destinations. A vehicle is modeled as consisting of two parts. The first part is the platform, which contains steering, engine, and wheels. The second part is a functional module that determines the purpose of the vehicle - either passenger transportation or freight delivery. The module can be exchanged at a depot or a special service depot, and thus change the purpose of the vehicle (see Figure \ref{fig:illustration}).In the proposed problem two different modules can be used, one exclusively for passenger transport (bus/ride-pooling operations) and the other for the exclusive delivery of freight items (freight truck operation). 

The contribution of this paper is threefold. First, the work addresses the research gap of sequentially integrated transportation systems for passenger and freight transportation on road networks. Second, the proposed problem formulation extends the existing research on PDP by separating the assignment of platforms and modules to requests, which allows the modelling of multiple changing processes at different locations throughout a single route. For these problem-specific model characteristics additional heuristics are implemented. Third, the paper studies the impact of multi-purpose vehicles through analyzing various large-scale scenarios.

The remainder of the paper is structured as follows. First, the relevant literature is reviewed in Section \ref{sec:literature_review}. Then the problem formulation and methodology of this study are described in Section \ref{sec:method}. The experimental design and the results are discussed in Section \ref{sec:experimental_design} and Section \ref{sec:results}, respectively. The paper closes with a critical assessment of the results, a conclusion, and an outlook for potential future studies in Section \ref{sec:conclusion}.
\section{Literature Review}\label{sec:literature_review}

The vehicle routing problem (VRP) has been widely studied for several decades and many variations in problem formulation and solution algorithms have been developed. Among the latest developments are four problem formulations relevant to the proposed formulation in this work: first, the heterogeneous vehicle routing problem as comprehensively reviewed by \citet{kocThirtyYearsHeterogeneous2016}; second, the swap body vehicle routing problem \citep{verologVeRoLogSolverChallenge2014}; third, the trailers and transshipment vehicle routing problem \citep{drexlApplicationsVehicleRouting2013}; and fourth, the share-a-ride problems. We elaborate on the relation between the problem addressed in this study and each of those in the subsequent paragraphs.

The heterogeneous routing problems describe logistic/routing operations with different types of demand and demand type specific vehicles to serve that demand. A common practical application of this routing problem is found in healthcare transport. In \cite{parraghIntroducingHeterogeneousUsers2011} and \cite{parraghModelsAlgorithmsHeterogeneous2012} the authors describe a DARP with heterogeneous demand and vehicles for the transportation of patients and disabled people. In their model vehicles may be equipped with staff seats, patient seats, stretchers and wheelchair places which in turn define the demand type and capacity for each vehicle. Another form of heterogeneous routing problems is considered by \cite{rekiekHandicappedPersonTransportation2006} and \cite{melachrinoudisDialarideProblemClient2007} where different vehicles in terms of capacity are utilized to serve a single type of demand. This problem is closely related to mixed vehicle routing problems, which do not consider pickup and drop off positions for each request. In comparison to these heterogeneous routing problems the model proposed in this work allows for an en-route change in vehicle configuration. In the previously studied heterogeneous routing problems, the vehicle configuration is decided upon before the depot departure and remains the same until the vehicle returns to the depot. Additionally, the number of configuration changes is not limited, allowing for several configurations for a given route. In \cite{qu_heterogeneous_2013} and \cite{tellez_fleet_2018} the authors present heterogeneous PDP and DARP with configurable vehicles, respectively. In these works vehicles can reconfigure their interior and, by that, change the capacity of the vehicle. The others propose a mixed-integer program which is solved using an ALNS in both papers. In \cite{qu_heterogeneous_2013} the authors analyze several scenarios and conclude that cost savings of 30\%-40\% can be achieved by changing the configuration of the vehicles. 

The Swap Body Vehicle Routing Problem (SB-VRP) was introduced as part of an operations research computation challenge and has been solved by several research teams, for example \citep{huberOrderMattersVariable2017, todosijevicGeneralVariableNeighborhood2017, toffoloStochasticLocalSearch2018}. In essence, the problem considers the routing of trucks, which can attach or remove trailers of a certain length. The nature of this problem is similar to the here proposed MP-PDP. However, in the SB-VRP the start and end depot for a trailer has to be one and the same, whereas in the MP-PDP a module can be loaded and dropped at any depot or service depot, hence embracing a more general and flexible functionality. Furthermore, no multi-depot functionality is implemented in the original problem formulation. Additionally, the SB-VRP can deal with two different types/sizes of trailers whereas the proposed work here can be easily extended to consider more vehicle types, e.g. passenger, freight, and waste transportation. Finally, the proposed vehicle routing problem extends the SB-VRP by adding additional constraints, such as maximum range per platform.

The third group of related vehicle routing problems are the trailer and transshipment problems \citep{drexlApplicationsVehicleRouting2013}. In addition to an adjusted objective function formulation, several additional constraints capturing the multi-depot considerations in the proposed MP-PDP formulation create a new problem variant. In truck-and-trailer routing problems only freight demand is typically considered (\cite{derigsTruckTrailerRouting2013} and \cite{parraghBranchandpriceAdaptiveLarge2017}).
Moreover, the types of trailer are limited to one; hence, only the addition or removal of trailers is considered. In contrast, the MP-PDP allows for the investigation of several demand type-specific modules, each of which having a different capacity.

\citet{liShareaRideProblemPeople2014} investigate if another mode of passenger transportation, namely private taxi rides, can be used for integrated urban transportation. The authors propose a reduced version of the Share-a-Ride problem and the Freight Insertion problem. The problem minimizes the additional operating costs of adding freight items in a set of planned taxi trips. The authors solve their proposed mixed-integer linear program (MILP) for static and dynamic demand scenarios. The numerical results are sensitive to the spatial distribution of the freight demand. The authors conclude that the integration of freight items into taxi services is a promising solution for urban areas. However, this new integrated mode should be complemented by a traditional truck service to guarantee the delivery of all packages. \cite{schroderIntegratedMultiagentUrban2017} present a multi-agent simulation model for passenger and freight transportation. They investigate the impacts of various policy measures, i.e. special vehicle tolls.

Since the computational complexity of vehicle routing and scheduling problems has proven to be NP-hard \citep{lenstraComplexityVehicleRouting22}, it is challenging to efficiently solve these problems for larger scenarios. Two general approaches are typically used to solve the VRP and its variations: (i) the utilization of exact analytical algorithms (e.g. Branch-and-Cut, Branch-and-Bound) and, (ii) the development of problem-specific heuristics or meta-heuristic algorithms (e.g. Simulated Annealing, Artificial Bee Colony, Genetic Algorithm or Large Neighborhood Search). In \cite{arslanCrowdsourcedDeliveryDynamic2016} a crowd-sourced delivery system for parcels and passengers is solved using an exact solution approach. An exact algorithm for the shared ride problem is presented by \cite{beirigoIntegratingPeopleFreight2018}, while \cite{ghilasPickupDeliveryProblem2016} solve the PDP with time windows and scheduled lines using CPLEX.

Due to the computational complexity of VRP problems, most researchers implement heuristic algorithms to solve large problems. \cite{chewGeneticAlgorithmBiobjective2013} develop a bi-objective genetic algorithm (GA) for the VRP and show its ability to reach improved objective values over previously published algorithms for Mandl's benchmark problems. \cite{alizadehforoutanGreenVehicleRouting2020} apply a similar genetic algorithm as well as simulated annealing (SA) algorithm to the green routing problem. The authors show that GA converges faster, whereas SA results in better solution robustness and qualities. In \cite{ropkeAdaptiveLargeNeighborhood2006} the adaptive large neighborhood search algorithm (ALNS) is presented. It is shown that the ALNS improves the best known solutions for VRP problems by around 50\%. Additional computational experiments indicate convergence robustness and its adaptability to various variations in problems. These results are confirmed in a later study by \cite{davidpisingerLargeNeighborhoodSearch2010} which shows that large-scale neighborhood search methods lead to fast and robust convergence for complex combinatorial problems. The authors propose variable local search algorithms and adaptive neighborhood definitions to further improve the computational efficiency of the algorithm. In the works of \cite{massonAdaptiveLargeNeighborhood2013}, \cite{ghilasAdaptiveLargeNeighborhood2016} and \cite{liAdaptiveLargeNeighborhood2016} the authors apply the ALNS to a variation of the VRPs. In their works, the authors solve the various VRPs by developing problem-specific heuristics.

Based on the review of the literature, we conclude that in recent years the integration and combination of multiple demand types in different transportation modes have been investigated. The benefits of which could be shown for simultaneous combination problems. Several authors proposed novel meta-/heuristic algorithms to solve complex VRP problems and showed their superiority in computation time and objective value over exact algorithms. Therefore, we have used a meta-heuristic optimization algorithm to solve the novel MP-PDP.

\section{Methodology}\label{sec:method}

The main characteristics of the pickup and delivery problem as proposed in this work are the multi-purpose vehicle concept and the optimal consolidation of passenger and freight items. In this section, the detailed problem formulation and solution algorithm are given.

\subsection{Proof of Concept}\label{sec:proof_of_concept}

\begin{figure*}[ht]
	\centering
	\begin{subfigure}[b]{0.49\textwidth}
		\centering
		\includegraphics[width=\textwidth]{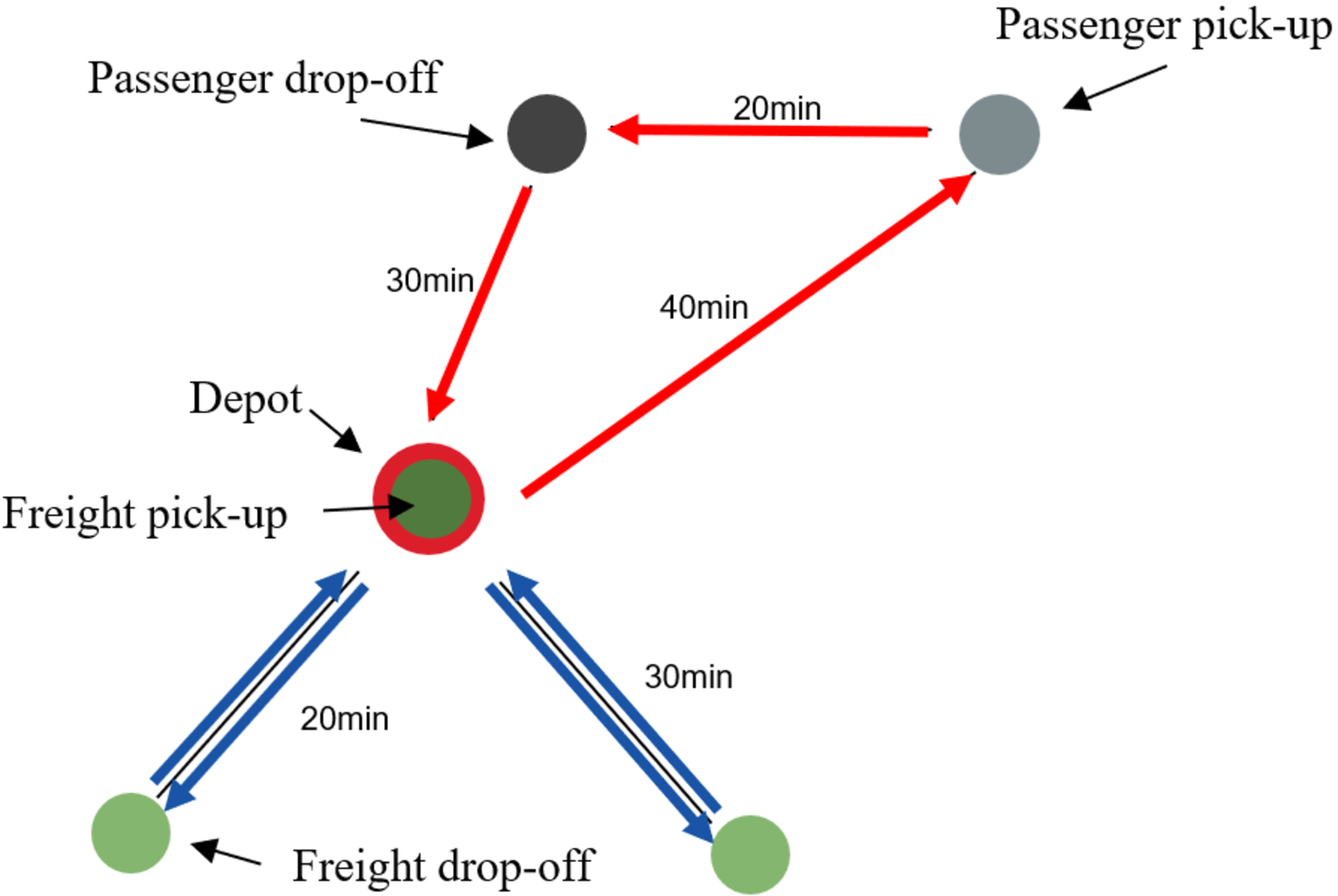}
		\caption[]%
		{Routing solution with conventional vehicles}
		\label{fig:conventional_proof}
	\end{subfigure}
	\hfill    
	\begin{subfigure}[b]{0.49\textwidth}
		\centering
		\includegraphics[width=\textwidth]{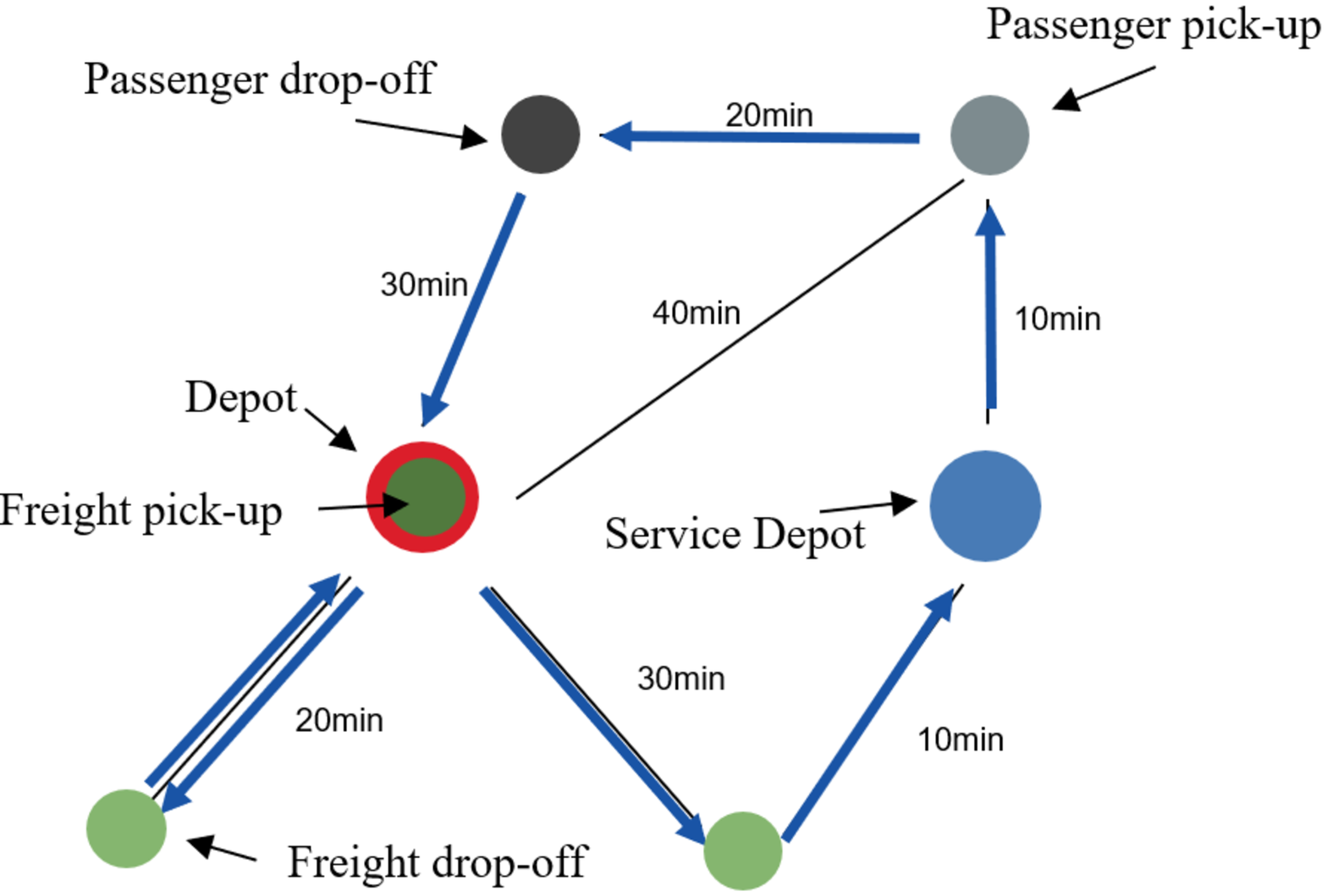}
		\caption[]%
		{Routing solution with multi-purpose vehicles}
		\label{fig:multi_proof}
	\end{subfigure}
	\caption{An illustration of conventional vehicle operations and multi-purpose vehicle operations}
	\label{fig:proof_of_concept}
\end{figure*}

In Figure \ref{fig:proof_of_concept} a simple routing example is given. This example showcases the theoretical benefits of multi-purpose vehicles over conventional vehicles.
Figure \ref{fig:conventional_proof} shows the conventional case. This solution is optimal and utilizes two vehicles (blue and red), one serving the freight requests, while the second vehicle serves the passenger requests. Therefore, this solution uses two platforms and two modules to serve the requests. In Figure \ref{fig:multi_proof} the same demand is served on a single route (blue). This solution utilizes the additional service depot to change the module and hence the purpose of that vehicle. Therefore, this solution uses one platform and two modules. The number of platforms could be reduced at the cost of exchanging modules once. Additionally, the total vehicle operation time is reduced from $20min + 20min + 30min + 30min + 40min +20min +30min = 190min$ to $20min + 20min + 30min + 10min + 10min +20min +30min = 140min$.

Finally, empty time can be reduced from $30min + 40min + 30min + 20min = 120min$ (either leaving or returning to the depot) to $30min + 10min + 10min + 20min = 70min$, thus producing a higher utilization of vehicle capacity, while serving all the demand at the same time as conventional vehicles. 


\subsection{Problem formulation}\label{sec:problem}

The problem formulation uses the nomenclature and parameter settings summarized in Table \ref{nomenclature}. The parameter values as stated in the Table \ref{nomenclature} are used in the experiments.

\begin{longtable}{c|l|c|c}
	\caption{Nomenclature and parameter values for the multi-purpose pickup and delivery problem} 
	\label{nomenclature} \\
Notation                      & Description        & Unit & Value  \\ 
\hline
$i,j$                           & Node index         & -    &  -      \\
$k  $                           & Platform index     & -    & -     \\
$m  $                           & Module index       & -    & -      \\
$u  $                           & Module usage index & -    & -      \\
$l  $                           & Depot group index  & -    & -      \\
\hline
$x_{i, j, k}$                &  Binary decision variable platform assignment                  & -    & -      \\
$y_{i, j, m}$                &  Binary decision variable module assignment                  & -    & -      \\
$s_{i, k}   $                &  Continuous decision variable platform arrival time                  & -    & -      \\
$e_{k, u}   $                & Binary decision variable module usage                   & -    & -      \\
$c_{i}      $                &   Continuous decision variable module capacity                 & -    & -      \\
\hline
$N       $                      &  \makecell{Set of all nodes in the graph (requests, depots and service depots) \\ $ND = N_r \cup N_d \cup N_{sd}$}                & -    & -      \\
$N_r     $                      &   Set of all request nodes $N = N^+ \cup N^-$                & -    & -      \\
$N^+_r     $                      &   Set of all requests pickup nodes                 & -    & -      \\
$N^-_r     $                      &   Set of all requests drop off nodes                 & -    & -      \\
$N^+_d     $                      &   Set of all origin nodes                  & -    & -      \\
$N^-_d     $                      &  Set of all destination nodes                   & -    & -      \\
$N_{sd}     $                      &    Set of all service depot nodes                & -    & -      \\
$M       $                      & Set of all available modules $M = M_p \cup M_f$                   & -    & -      \\
$M_p     $                     &    Set of all available passenger modules                &  -    & -      \\
$M_f     $                     &    Set of all available freight modules                &  -    & -      \\
$K       $                      &   Set of all available platforms                 &  -    & -      \\
$G_l     $                    & Set of depots in depot group $l \in N^+_d$         &  -    & -      \\
\hline
$a_i      $                   & Lower time window bound for node $i \in ND$        &   sec   & -      \\
$b_i      $                   & Upper time window bound for node $i \in ND$        &   sec     &      -  \\
$t_{i, j} $                 & Travel time including service time from node $i$ to $j$ for $ i,j \in ND$      &  sec    &    -    \\
$q_i      $                    & Demand for node $i \in ND$                   &  -    &-        \\
$w_{i, j} $                  &   Travel distance from node $i$ to $j$ for $ i,j \in ND$      &  m     &    -    \\
\hline
$h_r$                         &     Number of requests               &  -    &   -     \\
$h_d$                          &         Number of depots           &  -    &    -    \\
\hline
$v  $                           &     Travel speed               &  km/h    &  20      \\
$\kappa$         &   Max. number of platforms                 &  -    &     -   \\
$\mu_p$         &   Max. number of passenger modules                 &  -    &  -      \\
$\mu_f$         &   Max. number freight modules                 &  -    &    -    \\
$\eta$           &  Max. range per platform                  &  km    &   100     \\
$\vartheta$      &  Max. service depot visits per depot                  &   -   &    5    \\
$\gamma_p$      &  Max. capacity for passenger modules                  &  pas.    &    16    \\
$\gamma_f$      &  Max. capacity for freight modules                  &   items   &     16   \\
$\zeta$          & Large positive Number                   & -      &      -  \\
\hline
$\alpha_{tt}$ &   Cost of total vehicle travel time                 & EUR/h     &     6.9   \\
$\alpha_{ps}$ &  Cost of a platform                   & EUR     & 313.67      \\
$\alpha_{ms}$ & Cost of a module                   & EUR      &     156.84  \\
$\alpha_{td}$ & Cost of travel distance                   &  EUR/km    &    0.1    \\
$\alpha_{mc}$ & Cost of module change                   &  EUR    &  8.8       \\
$\alpha_{ud}$ & Cost of one unserved request                   &  EUR    & 470.52        
\end{longtable}

The MP-PDP is formulated as follows. Let $G = (N,A)$ be a directed graph where $N$ is a set of nodes and $A$ is a set of arcs. $ND = N_{d} \cup N_r \cup N_{sd}$ is the union of all depot nodes $N_d$, request nodes $N_r$, and the service depot nodes $N_{sd}$. The set of depot nodes is defined as $N_d = N^+_d \cup N^-_d$, where $N^+_d$ are the origin depot nodes and $N^-_d$ are the destination depot nodes. For each depot, there is one pair of nodes from $N_d$ in graph $G$. Similarly, is the set of request nodes defined as $N_r = N^+_r \cup N^-_r$, with the pickup nodes $N^+_r = N^+_{r,p} \cup N^+_{r,f}$ for passenger and freight requests, respectively, and the drop off nodes $N^-_r = N^-_{r,p} \cup N^-_{r,f}$ for passenger and freight requests, respectively. To simplify the mathematical problem formulation (compare Eqs. \eqref{obj} - \eqref{const:43}) is each pair of depot nodes duplicated by the maximum number of available platforms ($\kappa$). This ensures that at maximum all platforms can depart / arrive from / at the same physical depot but not at the same node in $G$. Similarly, the service depot nodes are duplicated by the maximum number of service depot visits ($\vartheta$), which is an input parameter determining how often the same service depot can be visited by a platform. This formulation implies that a maximum number of $\vartheta$ modules per type must be available at each physical service depot location to facilitate feasible operations. With $h_r = |N_r|/2$ and $h_d = \kappa \cdot |N_d|/2$ as the number of requests and the number of depots, respectively, $N$ can be defined as the set of nodes $N = \{1,...,h_d,h_d+1,...,h_d+h_r,h_d+h_r+1,...,h_d+2\cdot h_r,h_d+2\cdot h_r+1, ...,2\cdot h_d+2\cdot h_r,2\cdot h_d+2\cdot h_r+1,...,2\cdot h_d+2\cdot h_r+\vartheta\cdot n_{sd}\}$. The set of origin depot nodes ($N^+_d$) is then $\{h_1,...,h_d\}$, the set of destination depot nodes ($N^-_d$) is then $\{h_d+2\cdot h_r+1, ...,2\cdot h_d+2\cdot h_r\}$, the set of request pickup nodes ($N^+_r$) is then $\{h_d+1,...,h_d+h_r\}$, the set of request drop off nodes ($N^-_r$) is then $\{h_d+h_r+1,...,h_d+2\cdot h_r\}$, and the set of service depot nodes ($N_{sd}$) is then $\{2\cdot h_d+2\cdot h_r+1,...,2\cdot h_d+2\cdot h_r+\vartheta\cdot n_{sd}\}$. The corresponding pickup/drop off pair for one request can be identified using $\left(i, i+h_r\right) \, \forall i \in N^+_r$, while the corresponding origin/destination pair for one depot is $\left(i, i+h_d+2\cdot h_r\right) \, \forall i \in N^+_d$.

The set of available modules ($M$) is defined as $M = M_p \cup M_f$, where $M_p$ are the passenger modules and $M_f$ are the freight modules, with $M_p = \{1,...,\mu_p\}$ and $M_p = \{\mu_p + 1,...,\mu_p + \mu_f\}$, respectively. The set of all available platforms ($K$) is defined as $K = \{1,...,\kappa\}$, with $\kappa$ being the maximum number of platforms. 

Each node $i \in N$ has a demand $q_i$, service duration $o_i$, and request-specific time window ($a_i$, $b_i$). For depot nodes and service depot nodes the demand is set to zero. For pickup nodes, the demand represents the number of passengers or packages that should be picked up at that node. The sum of the demand for a given request pair is zero, i.e., the demand for drop off nodes is the negative value of the corresponding pickup node. The service duration for packages and passengers, and the respective time windows, are specified for each node individually. 

The set of all arcs $A$ is defined as $A = \{(i,j) \, \forall i,j \in N\}$ and hence spans a fully connected graph. The travel time $t_{i,j}$ is computed based as $t_{i,j} = d_{i,j}/v + o_i \, \forall i,j \in N$, where $v$ is the average travel speed in the network $G$. In our formulation, we assume that the platforms drive with a constant and known speed.

A homogeneous fleet of maximum $\kappa$ vehicle platforms is available, each using a(n) (exchangeable) module $m \in M$ with a type specific capacity. We assume that a fixed number of modules are available at each service depot at each point in time, defined by the parameter $\vartheta$. Vehicles experience a service time at a service depot and a regular depot which represents the time needed for any loading and/or unloading. The distance travelled by any individual platform cannot exceed the maximum range value $\eta$.

The decision variables for the problem are $x_{i,j,k}$, $y_{i,j,m}$, $e_{k,u}$, $s_{i,k}$, and $c_{i}$. $x_{i,j,k}$ is a binary variable which is 1 if the node pair $(i,j)$ is served by the platform $k$, and 0 otherwise. $y_{i,j,m}$ is a binary variable which is 1 if node pair $(i,j)$ is served by module $m$, and 0 otherwise. $e_{k,u}$ is a binary variable which is 1 if platform $k$ used modules $u$ and 0 otherwise. This variable is used to keep track of the total usage of the module. $s_{i,k}$ is a positive real number that indicates the arrival time of vehicle platform $k$ at node $i$. $c_{i}$ is a positive real number which indicates the current load at node $i$.

\newpage

The objective function parameters in the model are total vehicle travel time cost ($\alpha_{tt}$), the fixed cost of a platform ($\alpha_{ps}$), the fixed cost of a module ($\alpha_{ms}$), the cost of travel distance ($\alpha_{td}$), the fixed cost of module change ($\alpha_{mc}$), the fixed cost of one unserved request ($\alpha_{ud}$).

The objective function is composed of the minimization of six cost terms (see Equation \eqref{obj}).
\begin{align} \label{obj}
	\mathop {\min_{x,y,s,e,c} } {\kern 4pt}  
	& \alpha_{ps} \cdot B_1 + \alpha_{ms} \cdot B_2 + \alpha_{td} \cdot B_3 + \alpha_{tt} \cdot B_4 + \alpha_{mc} \cdot B_5 + \alpha_{ud} \cdot B_6
\end{align}

The first term in the objective function computes the total number of platforms used in a solution (see Equation \eqref{obj:1}). Note, that index $i$ is summed over the origin nodes $N^+_d$ of a scenario. Since an origin node is visited only by one platform, the total sum equals the platform fleet size.

\begin{equation}
\begin{aligned}
    B_1 &=\sum_{i \in N^+_d} \sum_{j \in N_r} \sum_{k \in K} x_{i, j, k} \label{obj:1}
\end{aligned}
\end{equation}

Equation \eqref{obj:2} computes the total number of modules used utilizing the module usage decision variable $e$ which is equal to $1$ if $u$ modules are used on platform $k$. The index $u$ is also used as an integer number with $\{ u \in \mathbb{Z} \, | \, 0 \leq u \leq |N| \}$.

\begin{equation}
\begin{aligned}
    B_2 &=\sum_{u \in N} \sum_{k \in K} u \cdot e_{k, u} \label{obj:2}
\end{aligned}
\end{equation}

Equation \eqref{obj:3} computes the total distance traveled by all platforms in the solution.

\begin{equation}
\begin{aligned}
    B_3 &=\sum_{i \in N} \sum_{j \in N} \sum_{k \in K} d_{i, j} \cdot x_{i, j, k} \label{obj:3}
\end{aligned}
\end{equation}

The fourth objective function term (see Equation \eqref{obj:4}) computes the total vehicle travel time time by summing the trip time for each platform, by subtracting the departure time from the origin depot from the arrival time of the corresponding destination depot. Note that for unused depots the subtraction equals zero, as is imposed by Equations \eqref{const:37} in the problem formulation.

\begin{equation}
\begin{aligned}
    B_4 &=\sum_{i \in N^+_d} \left( \sum_{k \in K} s_{i + 2 \cdot h_r + h_d, k}  - \sum_{k \in K} s_{i, k} \right) \label{obj:4}
\end{aligned}
\end{equation}

Equation \eqref{obj:5} computes the number of service depot visits, in combination with Equation \eqref{const:10} and \eqref{const:11} this equals the number of module changes.

\begin{equation}
\begin{aligned}
    B_5 &=\sum_{i \in N} \sum_{j \in N_{sd}} \sum_{k \in K} x_{i, j, k} \label{obj:5}
\end{aligned}
\end{equation}

The unserved demand is computed in Equation \eqref{obj:6} by subtracting the number of pickup nodes served from the total number of requests.

\begin{equation}
\begin{aligned}
    B_6 &= h_r - \sum_{i \in N^+_r} \sum_{j \in N} \sum_{k \in K} x_{i, j, k} \label{obj:6}
\end{aligned}
\end{equation}

\newpage

In addition to the objective function, the following constraints constitute the proposed MP-PDP.

\begingroup
\allowdisplaybreaks
\begin{align}
    \sum_{j \in N} \sum_{k \in K} x_{i, j, k}  &\leq 1  \qquad \forall i \in N, \label{const:1} \\
    \sum_{k \in K} x_{i, j, k} - \sum_{m \in M} y_{i, j, m} &= 0 \qquad \forall i \in N, j \in N, \label{const:2} \\
    \sum_{k \in K} x_{i, j, k} - \sum_{m \in M_p} y_{i, j, m} &= 0 \qquad \forall i \in N_{r,p}, j \in N, \label{const:3} \\
    \sum_{k \in K} x_{i, j, k} - \sum_{m \in M_f} y_{i, j, m} &= 0 \qquad \forall i \in N_{r,f}, j \in N, \label{const:4} \\
    \sum_{j \in N} x_{i, j, k} - \sum_{j \in N} x_{h_r + i, j, k} &= 0 \qquad \forall i \in N^+_r, k \in K, \label{const:5} \\
    \sum_{j \in N} y_{i, j, m} - \sum_{j \in N} y_{h_r + i, j, m} &= 0 \qquad \forall i \in N^+_r, m \in M, \label{const:6} \\
    \sum_{j \in N} x_{j, i, k} - \sum_{j \in N} x_{i, j, k} &= 0 \qquad \forall i \in N_r, k \in K, \label{const:7} \\
    \sum_{j \in N} y_{j, i, m} - \sum_{j \in N} y_{i, j, m} &= 0 \qquad \forall i \in N_r, m \in M, \label{const:8} \\
    \sum_{j \in N} x_{j, i, k} - \sum_{j \in N} x_{i, j, k} &= 0 \qquad \forall i \in N_{sd}, k \in K, \label{const:9} \\
    \sum_{j \in N} \sum_{m \in M_p} y_{j, i, m} - \sum_{j \in N} \sum_{m \in M_f} y_{i, j, m} &= 0 \qquad \forall i \in N_{sd}, \label{const:10} \\
    \sum_{j \in N} \sum_{m \in M_f} y_{j, i, m} - \sum_{j \in N} \sum_{m \in M_p} y_{i, j, m} &= 0 \qquad \forall i \in N_{sd}, \label{const:11} \\
    \sum_{j \in N} \sum_{k \in K} x_{i, j, k}  &\leq 1 \qquad \forall i \in N^+_d, \label{const:12} \\
    \sum_{j \in N} x_{i, j, k} - \sum_{j \in N} x_{j, \left(i + 2 \cdot h_r + h_d\right), k} &= 0 \qquad \forall i \in N^+_d, k \in K, \label{const:13} \\
    s_{i, k} + t_{i, j} - \zeta \cdot \left(1 - x_{i, j, k}\right) &\leq s_{j, k} \qquad \forall i \in N, j \in N, k \in K, \label{const:14} \\
    s_{i, k} + t_{i, h_r + i} - \zeta \cdot \left(1 - \sum_{j \in N} x_{i, j, k} \right) &\leq s_{h_r + i, k} \qquad \forall i \in N^+, k \in K, \label{const:15} \\
    s_{i, k} &\geq a_{i} \qquad \forall i \in N, k \in K, \label{const:16} \\
    s_{i, k} &\leq b_{i} \qquad \forall i \in N, k \in K, \label{const:17} \\
    \sum_{i \in N} \sum_{j \in N} w_{i, j} \cdot x_{i, j, k} &\leq \eta \qquad \forall k \in K, \label{const:18} \\
    c_{i} &\leq \gamma_p \qquad \forall i \in N_{r,p}, \label{const:19} \\
    c_{i} &\leq \gamma_f \qquad \forall i \in N_{r,f}, \label{const:20} \\
    c_{j} + \zeta \cdot \left(1 - y_{i, j, m} \right) &\geq c_{i} + q_{j} \qquad \forall i \in N, j \in N_r, m \in M, \label{const:21} \\
    c_{i} &= 0 \qquad \forall i \in \left(N^+_d \cup N_{sd}\right), \label{const:23} \\
    \sum_{u \in N} u \cdot e_{k, u} - \sum_{i \in \left(N^+_d \cup N_{sd}\right)} x_{i, j, k} &= 0 \qquad \forall j \in N, k \in K, \label{const:24} \\
    \sum_{i \in G_{l}} \sum_{j \in N} x_{i, j, k} &\leq 1  \qquad \forall l \in N^+_d, k \in K, \label{const:38} \\
    x_{i, i, k} &= 0 \qquad \forall i \in N, k \in K, \label{const:25} \\
    y_{i, i, m} &= 0 \qquad \forall i \in N, m \in M, \label{const:26} \\
    \sum_{j \in N} \sum_{k \in K} x_{i, j, k}  &= 0 \qquad \forall i \in N^-_d, \label{const:27} \\
    \sum_{j \in N} \sum_{m \in M} y_{i, j, m} &= 0 \qquad \forall i \in N^-_d, \label{const:28} \\
    x_{i, j, k} &= 0 \qquad \forall i \in N^+_r, j \in \left(N_{sd} \cup N^-_d\right), k \in K, \label{const:29} \\
    y_{i, j, m} &= 0 \qquad \forall i \in N^+_r, j \in \left(N_{sd} \cup N^-_d\right), m \in M, \label{const:30} \\
    x_{i, j, k} &= 0 \qquad \forall i \in \left(N_{sd} \cup N^+_d\right), j \in N^-_r, k \in K, \label{const:31} \\
    y_{i, j, m} &= 0 \qquad \forall i \in \left(N_{sd} \cup N^+_d\right), j \in N^-_r, m \in M, \label{const:32} \\
    x_{i, j, k} &= 0 \qquad \forall i \in N_{sd}, j \in N_{sd}, k \in K, \label{const:33} \\
    y_{i, j, m} &= 0 \qquad \forall i \in N_{sd}, j \in N_{sd}, m \in M, \label{const:34} \\
    x_{i, j, k} &= 0 \qquad \forall i \in N^+_d, j \in N_{sd}, k \in K, \label{const:35} \\
    y_{i, j, m} &= 0 \qquad \forall i \in N_{sd}, j \in N^-_d, m \in M, \label{const:36} \\
    s_{i, k} + \zeta \cdot \sum_{j \in N} x_{j, i, k} &\geq b_{\left(i - 2 \cdot h_r - h_d \right)} \qquad \forall i \in N^-_d, k \in K, \label{const:37} \\
    x_{i,j,k} &= 0 \quad or \quad 1 \qquad \forall i,j \in N, k \in K, \label{const:39} \\
    y_{i,j,m} &= 0 \quad or \quad 1 \qquad \forall i,j \in N, m \in M, \label{const:40} \\
    e_{k,u} &= 0 \quad or \quad 1 \qquad \forall u \in N, k \in K, \label{const:41} \\
    s_{i,k} & \geq 0 \qquad \forall i\in N, k \in K, \label{const:42} \\
    c_{i} & \geq 0 \qquad \forall i \in N \label{const:43}
\end{align}
\endgroup

The constraint in Equation \eqref{const:1} guarantees that each node is visited at most once. Using the less-than relation here and in equations \eqref{const:12} not every node has to be visited, hence allowing for unserved demand. In equation \eqref{const:2} the platform assignment is connected with the module assignment to guarantee that each served node pair $i,j$ is visited by a platform and a module.

In equations \eqref{const:3} and \eqref{const:4} the request type and module type are required to match, that is, passenger requests have to be served with passenger modules and freight requests with freight modules. Constraints \eqref{const:5} and \eqref{const:6} guarantee that the pickup node and its corresponding delivery node are served by the same platform and module, respectively. Similarly, equations \eqref{const:7} and \eqref{const:8} ensure that the same platform and module that enters a node also leaves that node (i.e., flow conservation constraints).

In constraints \eqref{const:9} - \eqref{const:11} the handling of platforms and modules in service depots is described. Equation \eqref{const:9} assures that the same platform entering a service depot node also leaves that service depot, whereas the module and its module type have to change as formulated in \eqref{const:10} and \eqref{const:11}. 

Each platform is required to start at an origin depot (compare Equation \eqref{const:12}) and end its trip at the corresponding destination depot (see Equation \eqref{const:13}). In Equation \eqref{const:14} the arrival times for consecutively visited nodes are constrained. Note that the travel time $t_{i,j}$ between a node pair $i,j$ includes the service time for node $i$. The large positive integer $\zeta$ is defined as 
$\max\left(b_i + t_{i,j} - a_j  \, \forall i,j \in N \right)$. 

In constraint \eqref{const:15} it is defined that the pickup node is served before the corresponding drop off node. Equations \eqref{const:16} and \eqref{const:17} guarantee that each arrival time is within its time window. In constraint \eqref{const:18} the length of each platform trip is constrained by the platform range. 

The set of constraints \eqref{const:19} - \eqref{const:23} guarantee that each module type is only filled up to its type-specific capacity (compare Equations \eqref{const:19} and \eqref{const:20}). Additionally, the capacity conservation is assured by equation \eqref{const:21}. Note that drop off nodes have the negative demand of its corresponding pickup node. In equation \eqref{const:23} the module loads are initialized to zero. 

Using equation \eqref{const:24} the number of used modules can be computed. The number of origin depot node $N^+_d$ or service depot node $N_{sd}$ visits by a platform $k$ is equal to the number of modules used. Note that 
$\{ u \in \mathbb{Z} \, | \, 0 \leq u \leq |N| \}$ is used as an index and integer variable. Due to duplication of depot nodes, the constraint \eqref{const:1} does not prevent the same platform from departing from the same physical depot twice. Therefore, additional constraints are required to ensure that the same platform $k$ does not depart from the same depot twice; in equation \eqref{const:38} this is formalized. $G_l$ contains a list of all duplicate nodes in the origin depot $l \in N^+_d$. 

In the constraints \eqref{const:25} - \eqref{const:37} several simplifications are formalized. Equations \eqref{const:25} and \eqref{const:26} prevent looping at any node. Equations \eqref{const:27} and \eqref{const:28} prevent any departure from any destination node for platforms and modules, respectively. Equations \eqref{const:29} and \eqref{const:30} prevent arriving at a service depot or destination node after visiting a pickup node for platforms and modules, respectively. Equations \eqref{const:31} and \eqref{const:32} prevent arriving at a drop off node after visiting a service depot or origin depot node for platforms and modules, respectively. Equations \eqref{const:33} and \eqref{const:34} prevent traveling from any service depot to any other service depot for platforms and modules, respectively. Equation \eqref{const:35} prevents traveling to any service depot after an origin depot node and equation \eqref{const:36} prevents traveling from any service depot to a destination depot node. In equation \eqref{const:37} the arrival times for unserved destination nodes are set to the upper time window bound of the corresponding origin depot node using $\zeta$ as a large positive number. By this, it is ensured that unserved depot nodes are not contributing to the trip duration term (compare Eq. \eqref{obj:4}). 

In the remaining equations \eqref{const:39} - \eqref{const:43} the domains of the decision variable for platform assignment ($x_{i,j,k}$), module assignment ($y_{i,j,m}$), module usage ($e_{k,u}$), arrival times ($s_i$) and capacity ($c_i$) are defined, respectively. 

\subsection{Adaptive Large Neighborhood Search}

For small and medium-sized instances, the problem can be solved using the optimization software CPLEX. For larger instances, a heuristic optimization algorithm is developed. We have adopted and implemented the adaptive large neighborhood search (ALNS) algorithm, originally proposed by \citep{ropkeAdaptiveLargeNeighborhood2006}, with problem-specific heuristics. The ALNS has been successfully applied to several PDPs (e.g. \cite{sacramentoAdaptiveLargeNeighborhood2019, massonAdaptiveLargeNeighborhood2013}) and has shown good performance for complex combinatorial problems, in general. The basic idea of the ALNS algorithm is to iteratively destroy and repair solutions. In Algorithm \ref{alg:1} the general outline of the implemented ALNS is given. 

\begin{algorithm} 
	\SetKwComment{Comment}{$\triangleright$\ }{}
	\KwData{passenger/freight demand, (service) depot positions, parameter settings}
	\KwResult{best solution ($x^*$) for the MD-PDP}
	create a feasible solution ($x$), set $x^* := x$\;
	\Repeat{maximum number of iterations, or objective variation threshold}{
		roulette wheel selection (see Equation \eqref{operator_propability}) for a destroy \& a repair operator using weights\;
		create a destroyed solution ($x_d$) using the chosen destroy operator on $x$\;
		create a candidate solution ($x'$) using the chosen repair operator on $x_d$\;
		\eIf{Objective($x'$) $<$ Objective($x^*$)}
		{
		    set $x^* := x'$\;
		    set $x := x'$\;
			set score for chosen destroy \& repair operator to $\sigma_1$\;
		}
		{
		    \eIf{$x'$ is accepted (Simulated Annealing)}
		    {
		        set $x := x'$\;
			    \eIf{Objective($x'$) $<$ Objective($x$)}
			    {
			        set score for chosen destroy \& repair operator to $\sigma_2$\;
			    }
			    {
			    	set score for chosen destroy \& repair operator to $\sigma_3$\;
			    }
		    }
		    {
		        set score for chosen destroy \& repair operator to $\sigma_4$\;
		    }
		}
		update weights using new operator scores (see Equation \eqref{weight_udpate})
	}
	return $x^*$
	\caption{An outline of the ALNS framework (see \cite{pisingerGeneralHeuristicVehicle2007})}
	\label{alg:1}
\end{algorithm}

In the first step, an initial feasible solution $x$ is created. This can be done using one of the repair operators, as all of them result in a feasible solution. The operator weights are all initialized to 1. In the next step, one destroy and one repair operator are selected using the roulette wheel selection process. The probability $p_{i,j}$ of choosing operator $i \in O$, with $O$ being the list of available operators, in iteration $j$ is calculated using Equation \eqref{operator_propability}. Here, $w_{i,j}$ is the weight for operator $i$ at iteration $j$. Note that there is a separate list ($O$) for destroy and repair operators, hence their probabilities are independent. The destroy and repair operators implemented are detailed below.

\begin{equation}
\begin{aligned}
    p_{i,j} = \frac{w_{i,j}}{\sum_{k \in O} w_{k,j}} \label{operator_propability}
\end{aligned}
\end{equation}

A candidate solution $x'$ is created by sequentially applying both operators. Note that after both operators have been applied, the candidate solution is always feasible. In the last step, the current solution $x$, global best solution $x^*$ and the heuristic scores are updated according to the objective value of $x'$ and if the candidate solution $x'$ is accepted. If a new overall best solution is found then the candidate solution is set as the new global best solution and the new current solution, and the score for the operators is set to $\sigma_1$. If the candidate solution is accepted using a simulated annealing decision process the current solution is set to the candidate solution. If the objective value of the accepted candidate solution is lower than the objective value of the current solution the operator scores are set to $\sigma_2$, else the score is set to $\sigma_3$. If the candidate solution $x'$ is not accepted, the best global solution and current solutions remain unchanged, and the operator scores are set to $\sigma_4$. At the end of each iteration $j$ the operator weights ($w_{i,j}$) for the chosen operators $i$ are updated using Equation \eqref{weight_udpate}. Here, the score ($s_{i,j}$) for the operator chosen $i$ at the current iteration $j$ is set to the corresponding $\sigma$ value as outlined in Algorithm \ref{alg:1}. The values of $\sigma$ and the computation of the updated weights is in line with the process described in \cite{ropkeAdaptiveLargeNeighborhood2006}.

\begin{equation}
\begin{aligned}
    w_{i,j+1} &= w_{i,j}\cdot \delta + \left( 1 - \delta \right) \cdot s_{i,j}   \label{weight_udpate}
\end{aligned}
\end{equation}

In Equation \eqref{weight_udpate} the operator decay parameter $\delta$ with $\{ \delta \in \mathbb{R} \, | \, 0 \leq \delta \leq 1 \}$, influences the speed at which the operator weights are adjusts to the scores. For low $\delta$ the adjustment rate is fast, while for high $\delta$ the adjustment rate is slow and the computed score values do not influence the operator's weights much. The algorithm terminates once the maximum number of iterations ($\lambda$) is reached or the change in the objective value - after a minimum number of iterations ($\lambda_{min}$) is performed - over an iteration span ($\omega$) is below a threshold ($\epsilon$). The termination criteria is computed using Equation \eqref{variation} at each iteration $i$.

\begin{equation}
\begin{aligned}
    \frac{\sum^{i-\omega}_{j = i - 2\cdot\omega}z_j}{\sum^{i}_{j = i - \omega} z_j} - 1 &\leq \epsilon \quad for, \lambda_{min} \leq i \leq \lambda, \label{variation}
\end{aligned}
\end{equation}

with $z_j$ being the objective value at each iteration.

\subsubsection{Heuristic Operators}

The implemented operators are in line with the conventionally used operators of \cite{ropkeAdaptiveLargeNeighborhood2006} and \cite{sacramentoAdaptiveLargeNeighborhood2019}. The interested reader is referred to these publications for more detailed information. In the following paragraphs only a brief description of each implemented operator and potential adjustments for the special nature of the MP-PDP are given.

\paragraph{Destroy Operators}

The general idea of the destroy operators is to help diversify the search process and thereby explore the solution space. In the implemented destroy heuristics depots (origin and corresponding destination node), service depots or requests (a pair of pickup and drop off nodes) are removed from a solution. Hence, except for service depots, no single nodes are removed from the solution but rather a node pair. 

\begin{itemize}
	\item Random removal: In this heuristic, a random selection of currently served requests (pickup and drop off node) are removed from the solution and considered as unserved.
	\item Module removal: In this heuristic a random selection of currently used modules is made. All requests served by any of the chosen modules, including potential service depot nodes, are removed from the solution and considered as unserved. 
	\item Platform removal: In this heuristic, one of the currently used platforms is chosen at random. All requests served by this platform, including potential service depot nodes, are removed from the solution and considered as unserved.
	\item Service depot removal: In this heuristic a random selection of currently used service depot nodes are removed from the solution. A service depot node can only be removed if that removal results in a feasible solution. Hence, in this heuristic redundant service depots such as two consecutive service depots, or service depot before/after a depot are removed. The removal of service depots including the removal of associated platforms is achieved with the module removal and platform removal heuristics. 
	\item Shaw removal: This removal heuristic was first proposed by \cite{shawNewLocalSearch1997a} and removes similar requests from the solution. In the work outlined here, the relatedness is computed using Equation \eqref{similarity}, using a distance term, travel time term, and request load term. To diversify the removal process, a determinism parameter is introduced; see detailed description in \cite{ropkeAdaptiveLargeNeighborhood2006}.
	\item Worst removal: In this heuristic, the requests with the highest cost contribution to the total cost are removed. The cost contribution for each served request is calculated by computing the objective function with and without this request.
\end{itemize}

\newpage

The relatedness ($R_{i,j}$) for two requests $i$ and $j$, as utilized in the \emph{Shaw removal} operator is computed with
\begin{equation}
\begin{aligned}
        R_{i,j} &= \phi \cdot \left(w_{a_i,a_j} + w_{b_i,b_j} \right) + \chi \cdot \left( |s_{a_i}-s_{a_j}| + |s_{b_i} - s_{b_j}| \right) + \psi \cdot \left( | q_i - q_j | \right), \label{similarity}
\end{aligned}
\end{equation}

where $\phi$, $\chi$, and $\psi$ are input parameter for the distance term, travel time term and the load term respectively. $a_i$ and $b_i$ are the pickup and drop off nodes for request $i$. $w_{i,j}$, $t_i$, and $q_i$ are the travel distance, arrival time, and demand for requests $i,j$, respectively.

The handling of service depot nodes introduces a few additional considerations that are novel to the proposed problem. In all destroy operators except for \emph{Module removal} and \emph{Platform removal}, a service depot node is only removed from the solution if that temporary solution remains feasible. Therefore, a removal would be feasible if a service depot is at the beginning or at the end of the route, or if two service depots are visited consecutively. In the\emph{Module removal} and \emph{Platform removal} operators, either an entire trip or a trip segment (i.e. depot to service depot, service depot to service depot, service depot to depot, or depot to depot) is removed. These operators strongly diversify the solution and guarantee that the number of platforms and modules is optimized. Additionally, these heuristics remove service depots from a solution, which is essential to further explore the solution space.  

In all destroy operators, except for the \emph{Platform removal}, the number of requests or modules which are removed is a random integer ($N$) following the discrete uniform distribution over the set 

\begin{equation}
\begin{aligned}
    N &= \{\iota, ..., min\left(n_{served}, n_r \cdot \xi \right)\}, \label{n_removal}
\end{aligned}
\end{equation}

with $\iota$ and $\xi$ being input parameters determining how many requests are removed. For high values of $\xi$ the maximum number of removed requests equals the number of served requests ($n_{served}$), while for low values of $\xi$ only a fraction of the total number of request ($n_r$) in the scenario are removed.

\paragraph{Repair Operators}

As a first step for all repair operators, the list of unserved requests is randomized. By this we avoid inserting the same requests in the same order multiple times. Next, we loop through the list of requests and try to insert each request individually following the principles described below. Additionally, we insert individual nodes (service depots) and node pairs (depots and requests) only in feasible positions in the temporary solution. This has consequences (1) at the end of each operator, we guarantee a feasible solution, and (2) certain solutions cannot be created immediately. An example of the second consequence is the situation when a request can only be inserted if a service depot is also added. In the current implementation of the algorithm, these cases are handled through the high number of iterations performed and the combination of intensification and diversification heuristics. 

At the end of each operator redundant service depots are removed. This applies in the following three cases: (1) a service depot is at the beginning of a route, (2) a service depot is at the end of a route, or (3) two service depots are visited consecutively.

\begin{itemize}
    \item First fit insert: In this repair heuristic, all unserved requests, depot nodes and service depot nodes are inserted into the first feasible location. In order to randomize this insertion process, the list of unserved nodes and the insertion routes are shuffled before each iteration.
    \item Inter route insert: In this repair heuristic all unserved requests, depot nodes and service depot nodes are inserted into the best feasible location of \emph{one} route. This route equals the route from which this request/node was previously removed. If the request/node has never been served, a random route is chosen. In order to randomize this insertion process, the list of unserved nodes is shuffled before each iteration. 
    \item Best insert: In this repair heuristic all unserved requests, depot nodes and service depot nodes are inserted into the best feasible location of \emph{all} routes. In order to randomise this insertion process the list of unserved nodes is shuffled before each iteration.
\end{itemize}

For \emph{Inter route insert} and \emph{Best insert} the best insertion position for each request/node is found using a two-step approach. In a first step, all feasible positions within a request route are computed. Then for each position, compute the objective value of the new route with the request/node inserted. The position with the smallest delta value is chosen as the final insertion position.

The re-/insertion of service depots in any of the three repair operators follows the following three principles. These principles are applicable to any route. (1) A service depot node can be inserted as the first node visited after the depot, (2) a service depot node can be inserted as the last node before ending the route, or (3) a service depot node can be inserted between any node pair, if the first node of that pair is a drop off node and the second node of that pair is a pickup node, and both these nodes are of different request type. The reconfiguration of platforms is achieved, since the list of requests and unserved (service) depots is randomized before each repair operator, and hence service depots can be inserted early or late in the repair process resulting in various new route configurations.

\section{Experimental Design}\label{sec:experimental_design}
To analyze the capabilities of the multi-purpose vehicle, several experiments are performed using the model described in the previous section.

\subsection{Scenario Definitions}
The scenarios are created by adjusting the spatial demand distribution, time window constraints, and number of service depots. In total 54 scenarios represent different vehicle use cases and operations. In the following paragraphs, the different scenario configurations are described. The analysis focuses on three main dimensions influencing the operation of the new vehicle concept: the temporal and spatial distributions of the passenger and freight demand, and the number of available service depots. Each scenario is solved using the parameters from Table \ref{nomenclature}. Scenarios without service depots model the conventional vehicle operations and function as the base case in the following analysis. The base demand for passengers and freight is in line with real data from the area of Stockholm, Sweden. For each scenario an ensemble run of 10 independent optimizations is performed, the average values of the ensemble runs are used for the analysis.

\begin{figure*}[ht]
	\centering
	\begin{subfigure}[b]{0.49\textwidth}
		\centering
		\includegraphics[width=\textwidth]{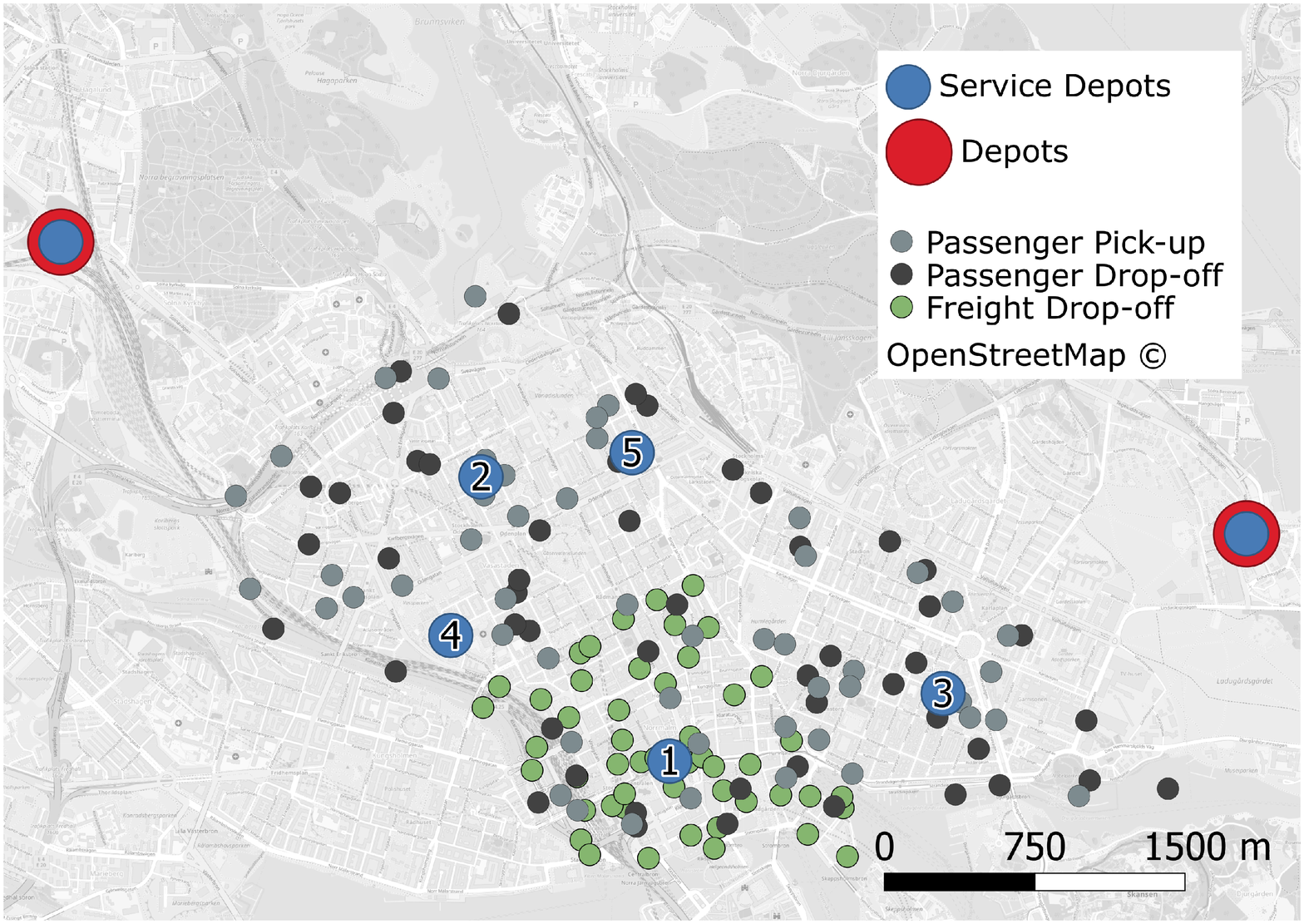}
		\caption[]%
		{Central representation (scenarios 1-18).}
		\label{fig:central}
	\end{subfigure}
	\hfill    
	\begin{subfigure}[b]{0.49\textwidth}
		\centering
		\includegraphics[width=\textwidth]{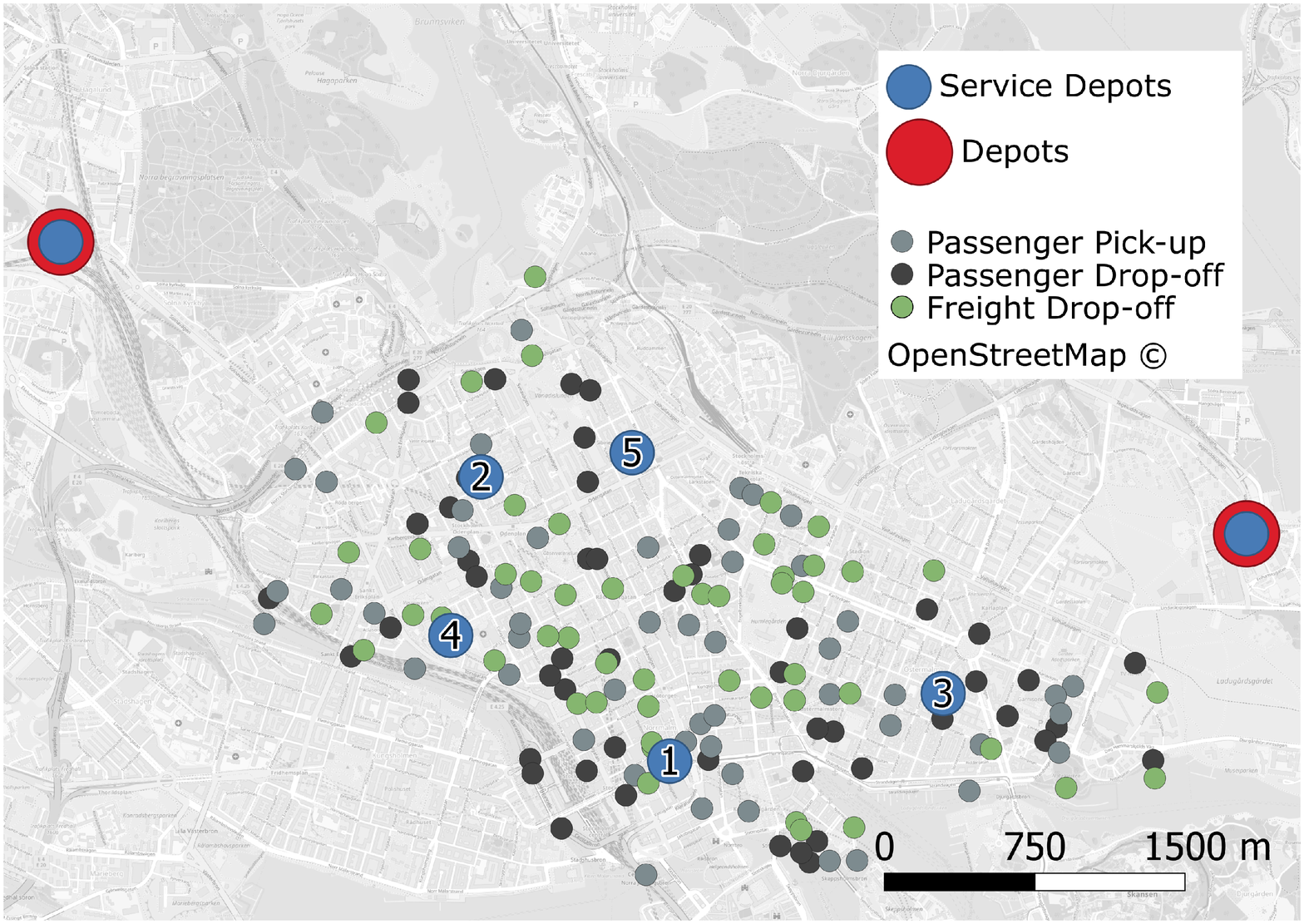}
		\caption[]%
		{Distributed representation (scenarios 19-36).}
		\label{fig:distributed}
	\end{subfigure}
	\hspace{0pt}    
	\begin{subfigure}[b]{0.49\textwidth}
		\centering
		\includegraphics[width=\textwidth]{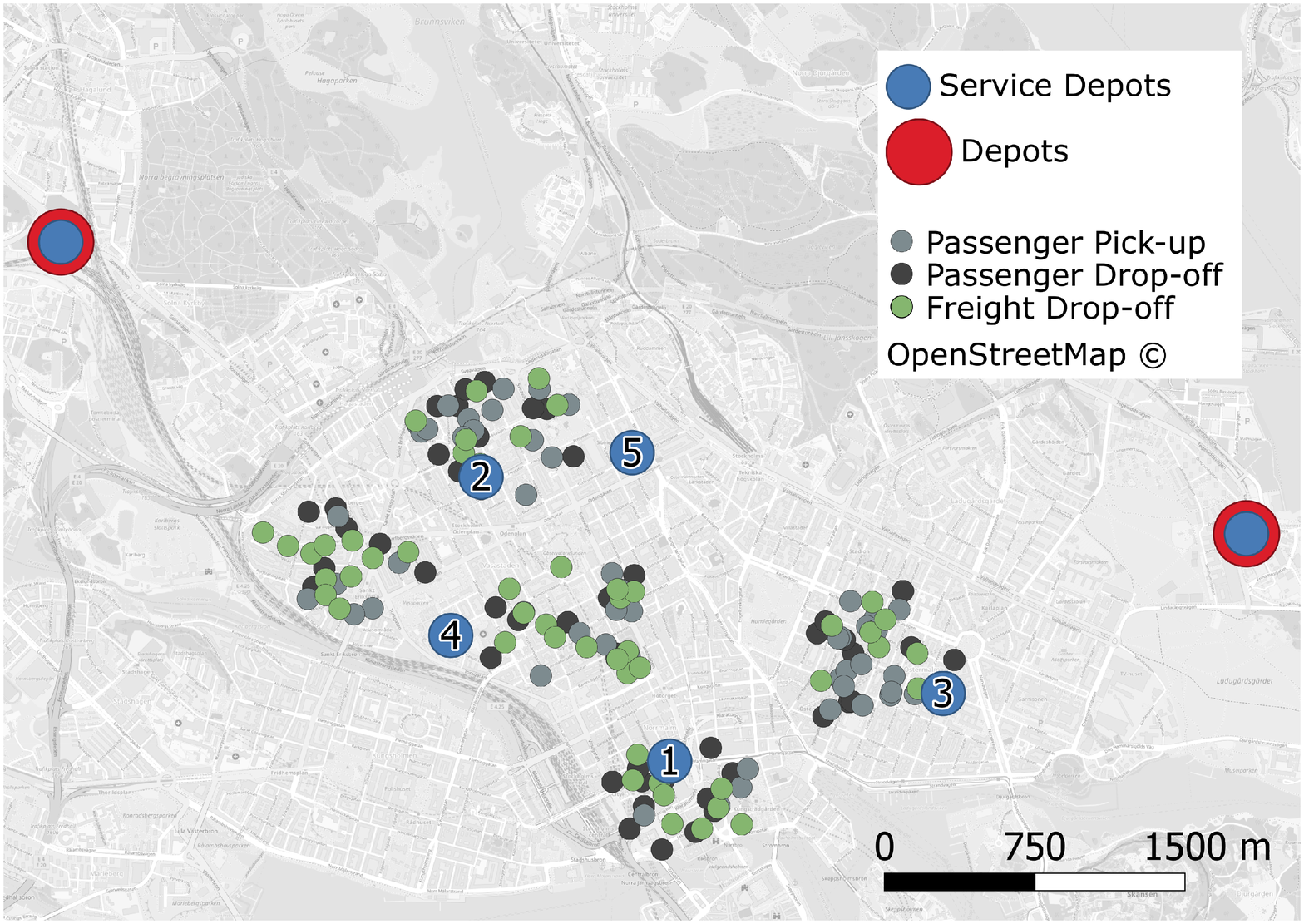}
		\caption[]%
		{Cluster representation (scenarios 37-54).}
		\label{fig:cluster}
	\end{subfigure}
	\caption{Spatial demand representation for the three different configurations.}
	\label{fig:scenarios}
\end{figure*}

\subsubsection{Spatial demand distribution}
Figures \ref{fig:central}--\ref{fig:cluster} represent the different spatial demand distributions analyzed in this paper. The demand distribution is distinguished in scenarios with central distribution (see scenarios 1-18 in Table \ref{tab:scenarios} and Figure \ref{fig:central}), scenarios with distributed distribution (see scenarios 19-36 in Table \ref{tab:scenarios} and Figure \ref{fig:distributed}) and scenarios with clustered distribution (see scenarios 37-54 in Table \ref{tab:scenarios} and Figure \ref{fig:cluster}). For all scenarios the demand consists of 50 passengers and 50 freight requests. The locations of the depots (red) correspond to those of freight distribution centers in Stockholm. The locations of service depots (blue) are chosen at strategic points in the service area to best illustrate the transport concept. The locations of depots and service depots remain unchanged among the different scenarios.

In scenarios 1-18 as described in Table \ref{tab:scenarios} (see Figure \ref{fig:central}) the passenger requests (pickup and drop off locations) are evenly distributed over the entire service area. The freight requests are focused on the central area in Stockholm. This central demand pattern represents urban deliveries operations that selectively serve e.g. stores, shops, and grocery stores only in a centralized area.

The spatial distribution of scenarios 19-36 from Table \ref{tab:scenarios} (see Figure \ref{fig:distributed}) represents city-wide deliveries spread evenly over the entire area of interest and typical passenger movement patterns throughout the city. In this scenario, the passenger and freight requests are spread out evenly over the entire service area. The passenger pickup and drop off points are different to Figure \ref{fig:central} but share the same generation process.

In the third variation in demand distribution (Figure \ref{fig:cluster} and scenarios 37-54 in Table \ref{tab:scenarios}), the passenger and freight requests form spatial clusters. A request can however be transported between clusters, meaning pickup and drop off locations do not need to be in the same cluster. An operator serving only a certain area of a city or having dedicated areas of interest would face a similar demand distribution as is represented in this scenario.

\subsubsection{Time window settings}
Next to the spatial demand variation, each scenario also has different time window definitions (compare Table \ref{time_window}). In total three different time window settings are implemented. The first setting represents a scenario where no time window constraints are present so that each request is not subject to a specific time. This setting is chosen to represent a base scenario for the optimization. Due to the removal of time windows the solution space is significantly increased and therefore finding a robust solution is more challenging for the algorithm. At the same time creating feasible solutions is simplified and the general functionality of the algorithm can be better shown. 

The second time window setting represents loose peak time window constraints, meaning that each request has a dedicated time for pickup and drop off but large deviations are permissible, e.g. several hours. The distribution of the requests follows a peak distribution, where passenger demand is highest during the morning and afternoon hours and freight demand is highest during mid-day.

In the last time window setting, the constraints are tightly defined. All requests have an arrival time window of \SI{20}{min}, while the departure time window is set to be unrestrictive. Similar to the second setting, the requests follow a peak distribution, i.e. passenger requests peak during the morning and afternoon and freight request peak during mid-day.
\begin{table}
	\caption{Time window definitions} 
	\label{time_window}
\centering
\begin{tabular}{c||c||cc||c}
Start Time [sec] & End Time [sec] & \multicolumn{2}{c||}{Type} & Count  \\ 
\hline\hline
No time windows& & & \\
1                & 86400          & ~drop off & ~freight     & 50     \\
1                & 86400          & ~drop off & ~passenger   & 50     \\
1                & 86400          & ~pickup  & ~freight     & 50     \\
1                & 86400          & ~pickup  & ~passenger   & 50    \\\hline
Peak time windows& & & \\
28800 & 43200 & ~drop off & ~passenger & 25  \\
28800 & 43200 & ~pickup  & ~passenger & 25  \\
36000 & 57600 & ~drop off & ~freight   & 50  \\
36000 & 57600 & ~pickup  & ~freight   & 50  \\
50400 & 64800 & ~drop off & ~passenger & 25  \\
50400 & 64800 & ~pickup  & ~passenger & 25 \\ \hline
Tight time windows& & & \\
1     & 86400 & ~pickup  & ~freight   & 50  \\
1     & 86400 & ~pickup  & ~passenger & 50  \\
24600 & 25800 & ~drop off & ~passenger & 2   \\
25800 & 27000 & ~drop off & ~passenger & 3   \\
27000 & 28200 & ~drop off & ~passenger & 5   \\
28200 & 29400 & ~drop off & ~passenger & 6   \\
29400 & 30600 & ~drop off & ~passenger & 5   \\
30600 & 31800 & ~drop off & ~passenger & 3   \\
31800 & 33000 & ~drop off & ~passenger & 1   \\
33000 & 34200 & ~drop off & ~passenger & 1   \\
34200 & 35400 & ~drop off & ~freight   & 1   \\
35400 & 36600 & ~drop off & ~freight   & 1   \\
36600 & 37800 & ~drop off & ~freight   & 1   \\
37800 & 39000 & ~drop off & ~freight   & 1   \\
39000 & 40200 & ~drop off & ~freight   & 2   \\
40200 & 41400 & ~drop off & ~freight   & 4   \\
41400 & 42600 & ~drop off & ~freight   & 4   \\
42600 & 43800 & ~drop off & ~freight   & 5   \\
43800 & 45000 & ~drop off & ~freight   & 6   \\
45000 & 46200 & ~drop off & ~freight   & 6   \\
46200 & 47400 & ~drop off & ~freight   & 5   \\
47400 & 48600 & ~drop off & ~freight   & 4   \\
48600 & 49800 & ~drop off & ~freight   & 4   \\
49800 & 51000 & ~drop off & ~freight   & 2   \\
51000 & 52200 & ~drop off & ~freight   & 1   \\
52200 & 53400 & ~drop off & ~freight   & 1   \\
53400 & 54600 & ~drop off & ~freight   & 1   \\
54600 & 55800 & ~drop off & ~freight   & 1   \\
55800 & 57000 & ~drop off & ~passenger & 1   \\
57000 & 58200 & ~drop off & ~passenger & 2   \\
58200 & 59400 & ~drop off & ~passenger & 3   \\
59400 & 60600 & ~drop off & ~passenger & 5   \\
60600 & 61800 & ~drop off & ~passenger & 6   \\
61800 & 63000 & ~drop off & ~passenger & 5   \\
63000 & 64200 & ~drop off & ~passenger & 2  
\end{tabular}
\end{table}

\subsubsection{Variation in available service depots}
The third dimension of the scenario definition is the number of available service depots. Each spatial and temporal configuration is combined with six different service depot variations, increasing from no available service depots to five available service depots in the city center (compare Figure 3). In all scenarios the same service depot is added first, second, etc. In the scenarios with at least one service depot in the city center the regular depots can also be used as service depots. The scenarios without available service depots represent the operations with conventional vehicles. 

In Table \ref{tab:scenarios} an overview of the configurations of all 54 scenarios is given.

\begin{table}[ht]
	\centering
	\caption{Scenario overview}
	\begin{tabular}{c|c|c|c}
		Scenario & Spatial distribution & Time window variation &        Nb. of service depots                   \\ \hline
		  1-6    &       Central        &      noTW     & (0, 1, 2, 3, 4, 5)  \\
		  7-12   &       Central        &  peak  & (0, 1, 2, 3, 4, 5)  \\
		 13-18   &       Central        &  tight  & (0, 1, 2, 3, 4, 5)  \\
		 19-24   &     Distributed      &      noTW     & (0, 1, 2, 3, 4, 5)  \\
		 25-30   &     Distributed      &  peak  & (0, 1, 2, 3, 4, 5)  \\
		 31-36   &     Distributed      &  tight  & (0, 1, 2, 3, 4, 5)  \\
		 37-42   &       Clustered        &      noTW     & (0, 1, 2, 3, 4, 5)  \\
		 43-48   &       Clustered        &  peak  & (0, 1, 2, 3, 4, 5)  \\
		 49-54   &       Clustered        &  tight  & (0, 1, 2, 3, 4, 5) 
	\end{tabular}
	\label{tab:scenarios}
\end{table}

\subsection{Parameter settings}

The parameters used for the experiments are determined based on parameter tuning procedures, reported literature values or general reasoning. If not stated otherwise the parameters in Table \ref{nomenclature} and Table \ref{alns_parameter} are used for the remainder of this paper. 

\begin{table}
	\caption{ALNS parameter} 
	\label{alns_parameter}
    \centering
    \begin{tabular}{l|l|l}
        Notation                        & Description  & Value  \\ 
        \hline
    $ \sigma_1 $        &    Operator weight for new best overall solution         &     7        \\
    $ \sigma_2 $        &     Operator weight for better candidate solution        &      2       \\
    $ \sigma_3 $        &   Operator weight for accepted solution        &       9      \\
    $ \sigma_4 $        &  Operator weight for rejected solution           &      1       \\
    $ \delta $           &      Operator decay per iteration       &       0.1       \\
    $ \lambda $          &    Iterations         &     10000        \\
    $ \lambda_{min}$ &      Min. iterations       &     2000         \\
    $ \omega $           &    Iterations look-back         &    1000          \\
    $ \epsilon $         &   Objective improvement threshold          &     0.01         \\
    $ T_{start} $                   &     Start temperature simulated annealing             &   100     \\
    $ T_{end}   $                   &     End temperature simulated annealing              &     0.0001   \\
    $ \nu $              &     Step size simulated annealing        &   0.9999           \\
    $ \phi $             &    Relatedness parameter for distance         &    9          \\
    $ \chi $             &     Relatedness parameter for time        &         3     \\
    $ \psi $             &  Relatedness parameter for load           &      2        \\
    $ \rho $             & Randomising selection parameter for shaw removal requests            &     6         \\
    $ \rho_{worst}$  &    Randomising selection parameter for worst removal requests         &     3         \\
    $ \xi $              &    Max. request removal factor         &     0.32         \\
    $ \iota$   &    Min. request removals         &      1       
    \end{tabular}
\end{table}

\subsection{Model parameters}

The platform and operational parameters, i.e. maximum trip range, and passenger/freight capacity are based on vehicle design parameters published by manufacturers or project reports. The travel time costs are based on a value-of-time of \SI{6.9}{EUR/h} as mentioned in \cite{borjessonExperiencesSwedishValue2014}. 

In \cite{militaoOptimalFleetSize2021} the authors present a linear estimation model based on vehicle capacity to estimate the distance- and time-based operational costs for electric vehicles. The model is build using data from Munich, Germany, and Santiago, Chile. For the distance related cost their model is $\alpha_{td} = 0.003599 \cdot \gamma_k + 0.04162$. For the time related costs their model equals $\alpha_{ps} = \left(0.1753 \cdot \gamma_k + 16.8\right) \cdot \frac{2}{3} \cdot 24$, and $\alpha_{ms} = \left(0.1753 \cdot \gamma_k + 16.8\right) \cdot \frac{1}{3} \cdot 24$, respectively, using the capacity $\gamma_k$ per type $k \in K$. The parameters $\alpha_{ps}$ and $\alpha_{ms}$ are scaled with $2/3$ and $1/3$, respectively, to account for cheaper module operation costs. The second factor, 24, represents the duration of the planning period in hours.


The module change cost is estimated based on an average salary of \SI{17.6}{EUR/h} for a bus driver in Sweden and the assumption that one person would need approximately \SI{30}{min} to switch the module. 

The parameter settings for unserved demand are chosen to match the fleet size related costs for one platform and one module. This is done so that a single vehicle is not used to serve a single request. 

There are no differences in model parameters between conventional operations and multi-purpose operations. This is done to isolate the impact of the new vehicle technology on the results.

In preparation for the final results, two experiments are performed. Both having the goal to overcome the challenges imposed by the numerical differences in model parameters and to achieve robust optimization of the vehicle fleet and vehicle routes. The first experiment uses a hierarchical objective function. The top-level problem optimizes the number of platforms and modules, while the low-level problem optimizes the vehicles routes using the results from the top-level problem. The objective function in this approach is adjusted to only represent the relevant cost terms for each subproblem. In the second experiment, the vehicle fleet and routes are optimized simultaneously using all cost terms in the objective function. When comparing the results of both experiments for all scenarios, no significant differences in the solutions could be seen. Hence, the results reported in Section \ref{sec:results} are achieved using the simultaneous optimization approach of the second experiment. 

Additionally, we employ a two-step approach. First, the scenarios with conventional vehicle operations are solved using the ALNS as described in Section \ref{sec:method}. Second, the best solutions of the first step are used as initial solutions for the scenarios with multi-purpose vehicle operations. Since all solutions with conventional operations are also feasible for the multi-purpose scenarios this improves the optimization process for the multi-purpose vehicles by initiating their optimization with high quality solutions. 

\subsection{Optimization parameter}

Following the principle of \cite{ropkeAdaptiveLargeNeighborhood2006} an iterative grid search parameter tuning approach has been performed. For each parameter in Table \ref{alns_parameter} a discrete set of 10 possible parameter values each has been manually specified. Additionally, 36 scenarios have been created synthetically with varying spatial/temporal demand distribution and number of (service) depots to have a diverse set of scenarios. Each scenario has a total of 40 requests, which represents a compromise between fast computation times and scenario complexity. The final parameter values are reported in Table \ref{alns_parameter}. These parameters correspond to the parameter values with the lowest average optimally gap over all scenarios compared to the best available solution.

\section{Results}\label{sec:results}

The result section is divided into two parts. In the first part, the ALNS and CPLEX solutions are compared for small and medium-sized problems to showcase the functionality and general usability of the ALNS algorithm. The second part presents and discusses the results obtained for large scenarios, focusing on the changes induced by the introduction of the multi-purpose vehicle technology.
 
\subsection{Validation of ALNS}

Small and medium-sized problems (see Table \ref{validation}) are used to validate the performance of the ALNS. The problems are created synthetically with varied number of requests, number of depots, and varied time window definitions. As shown in Table \ref{validation}, this variation leads to several problems for each value of the number of requests. In total 118 problems with varied complexity were created. 

To validate the ALNS performance, each problem was solved using CPLEX and ALNS. For the ALNS approach, each problem was solved with an ensemble run of 10. The numerical values reported for the ALNS for each request group are based on the average values of each ensemble run. The maximum computation time for CPLEX was set to \SI{1}{h} and the accepted optimization gap was set to 0.4\%.

As can be seen in Table \ref{validation}, the ALNS algorithm is capable of finding the optimal solution for all problems. The objective deviation is computed as the percentage difference between the best CPLEX solution and the average of all ensemble ALNS solution for a problem. The reported value in Table \ref{validation} is the objective deviation averaged over all problems with the same number of requests. With the CPLEX termination conditions six of the problems could not be solved to optimality within the permitted time limit, meaning that the ALNS is in all cases capable to match the performance of CPLEX. In fact, in 116 problems the ALNS always finds the exact optimal solution or better solutions than CPLEX due to exceeded time limit.

Looking at the average computation time for CPLEX and ALNS a clear trend can be seen. For larger problems starting from approximately 20 requests the ALNS computation time is shorter than that of CPLEX.  The reported computation time for ALNS is the time until the algorithm terminates, e.g., over 2000 iterations. However, for most problems, the final best solution is found after only a couple of hundred iterations, highlighting the rapid convergence of ALNS. It should be noted that the ALNS computation time is not linear with the number of requests.

\begin{table}
\caption{Comparison of CPLEX and ALNS}
\label{validation}
\centering
\begin{tabular}{c|c|c|c|c}
Nb. of requests & Nb. of problems & Avg. Obj. deviation & Avg. time CPLEX [sec] & Avg. time ALNS [sec]  \\ 
\hline
4           & 12   & -4.11E-06                 & 0.21                                 & 2.52                                 \\
8       &     12   & -2.13E-09                 & 0.48                                 & 2.81                                 \\
12      &   12     & -2.65E-06                 & 1.26                                 & 2.63                                 \\
16       &    12   & 8.03E-10                  & 2.62                                 & 3.11                                 \\
20        & 12     & -7.60E-07                 & 11.07                                & 4.43                                 \\
24         &  12   & 3.08E-08                  & 44.83                                & 4.35                                 \\
28      & 10       & 2.94E-08                  & 422.70                               & 5.04                                 \\
30      &   24     & 6.62E-03                  & 427.90                               & 5.16                                 \\
40      &   12     & 4.38E-08                  & 1707.35                              & 6.16                                
\end{tabular}
\end{table}

\subsection{Scenario analysis}
In this section scenario-specific results with a focus on the spatial, temporal, and service depot variation are discussed. The average computation time for all scenarios is \SI{1020}{sec}. The final best solution was found on average after 1200 iterations, while on average 2020 iterations were computed for each scenario.

\subsubsection{Impact of spatial distribution and time window constraints}

\begin{figure*}[ht]
	\centering
	\includegraphics[width=0.9\textwidth]{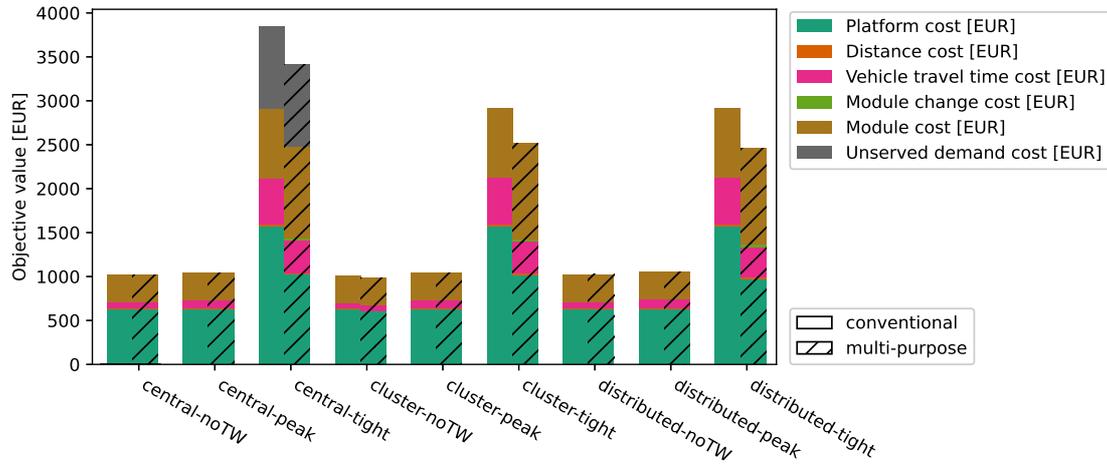}
	\caption{Average objective values for different spatial distributions and time window definitions.}
	\label{fig:Objective_aggregated}
\end{figure*}

The figures in this section show aggregated values. The numerical values are average values for each aggregation. In Figure \ref{fig:Objective_aggregated} the objective value for the different definitions of the time window and the spatial groups is shown. The bars without pattern correspond to scenarios without service depots and therefore represent conventional operations, while the multi-purpose operations have between 1 and 5 additional service depots and are shown as bars with pattern. The different colors indicate the value of each cost term in the objective function.

\begin{figure*}[ht]
	\centering
	\includegraphics[width=0.75\textwidth]{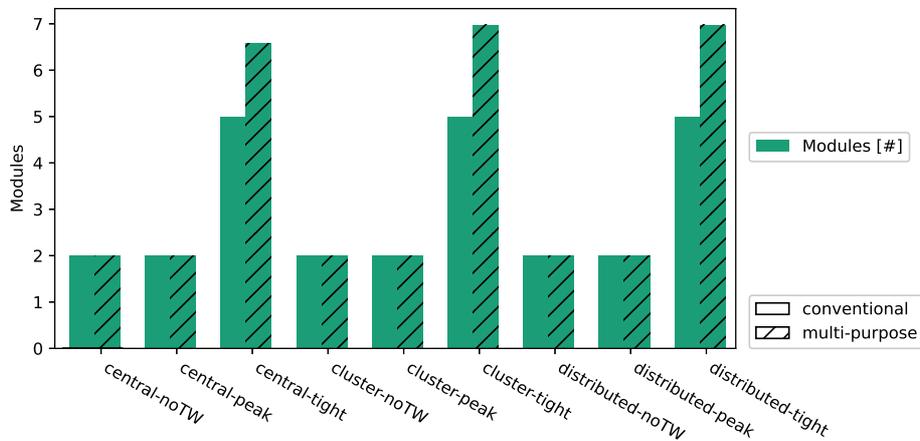}
	\caption{Average number of modules utilized for different spatial distributions and time window definitions.}
	\label{fig:Modules_aggregated}
\end{figure*}

\begin{figure*}[ht]
	\centering
	\includegraphics[width=0.75\textwidth]{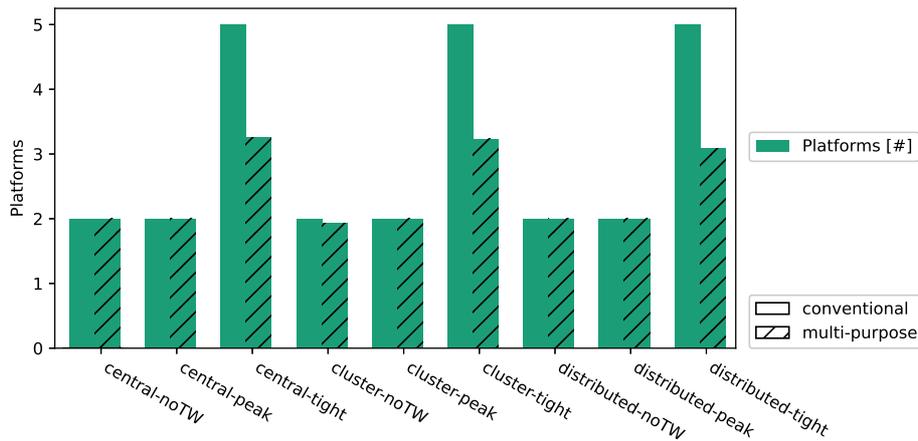}
	\caption{Average number of platforms utilized for different spatial distributions and time window definitions.}
	\label{fig:Platforms_aggregated}
\end{figure*}

It can be seen that for tight time window definitions the objective values are significantly higher than for peak and no time window definitions. This is because more platforms and more modules are needed to serve all requests. In addition to that, there is no significant difference in the objective value between operations with no time window and peak time window definitions. The objective value reduction are 11.3\%, 13.4\%, and 15.4\% for central, cluster and distributed scenarios, respectively. Meaning that multi-purpose operations are most beneficial for clustered and distributed scenarios. This reduction for tight time window definition mainly stems from a reduction of used platforms of approximately 35\%, which outweighs the higher costs for the modules required. In the conventional case additional vehicles have to be used to serve requests which cannot be served within the tight-time window, whereas the flexibility of the multi-purpose vehicles allows for a strategic change of the module just so that more requests can be served with the same platform. In peak time window scenarios this advantage is not given since most requests can be served in an optimal way with the same amount of platforms. This principle can clearly be seen on Figure \ref{fig:Modules_aggregated} and Figure \ref{fig:Platforms_aggregated}. Additionally, the total vehicle travel time costs and distance costs are reduced by approximately 33\% and 16\%, respectively, when deploying multi-purpose vehicles. The significantly higher objective values for centralized scenarios with tight time windows compared to other tight time window scenarios can be explained with the unserved demand. In all these scenarios, 2 requests could not be served. This is independent of the vehicle technology used. 

\begin{figure*}[ht]
	\centering
	\includegraphics[width=0.8\textwidth]{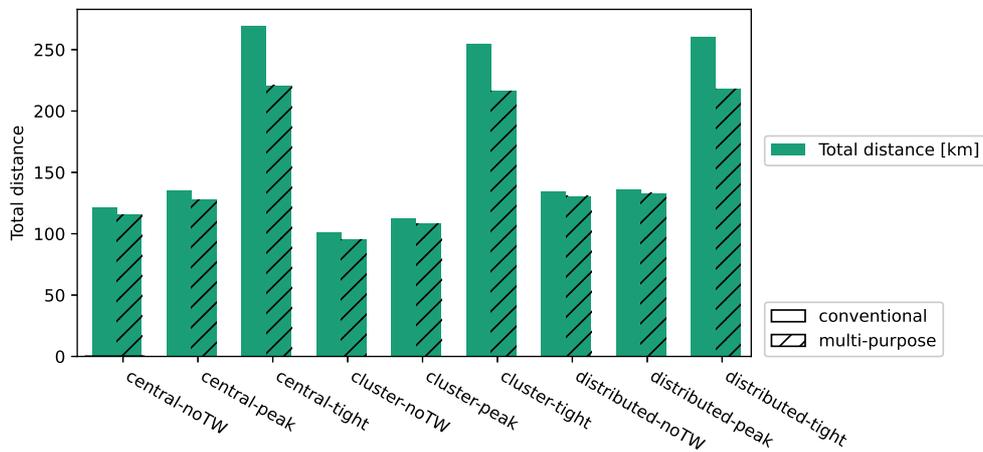}
	\caption{Total distance traveled by all platform in a solution averaged for different spatial distributions and time window definitions.}
	\label{fig:Total_distance_aggregated}
\end{figure*}

\begin{figure*}[ht]
	\centering
	\includegraphics[width=0.85\textwidth]{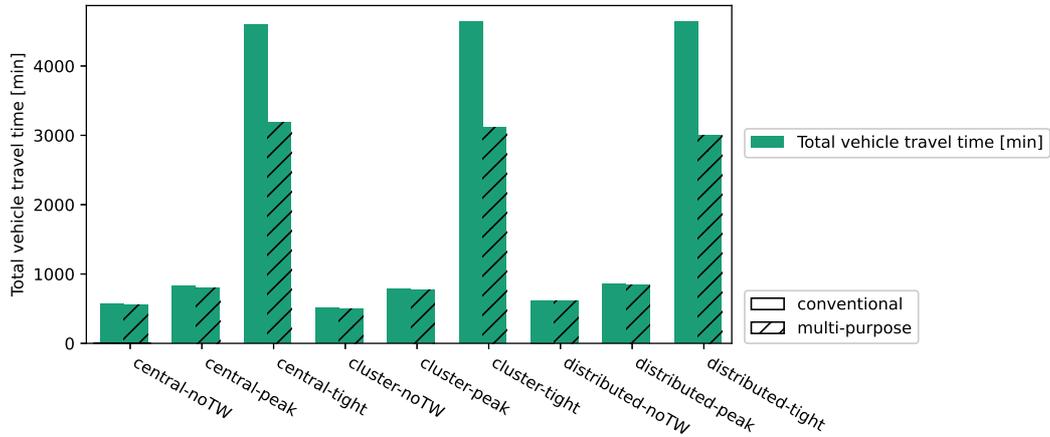}
	\caption{Trip duration by all platforms in a solution averaged for different spatial distributions and time window definitions.}
	\label{fig:Trip_duration_aggregated}
\end{figure*}

When analyzing the total distance traveled by all platforms in a solution (see Figure \ref{fig:Total_distance_aggregated}) and the total travel time cost (see Figure \ref{fig:Trip_duration_aggregated}), the following general observations can be made. First, for peak and no time window scenarios the values for total distance and trip duration are significantly lower than for the tight time window scenarios. This is as expected since the order in which requests can be served is less constraint, hence the ALNS algorithm can generate more direct and therefore shorter platform routes, which in turn also results in shorter trip durations. This observation is independent of the vehicle technology. Interestingly, for tight time window scenarios, the trip durations for multi-purpose vehicle operations are significantly shorter compared to conventional operations. This can be explained by the fact that the multi-purpose vehicles utilize the centrally located service depots to start a new trip, which directly translates to shorter trip times. This also results in a reduction in the distance traveled (see Figure \ref{fig:Total_distance_aggregated}). The explanation for this is two fold. First, to achieve shorter trip duration and trip distances the order in which requests are served is different between conventional and multi-purpose operations. This new request order can introduce shorter travel distances. Second, due to the removal of additional trips to and from the main depots shorter travel distances and trip times can be recorded. A second observation can be made by comparing the spatial distributions. For clustered scenarios without and peak time windows, the total distance traveled by the platforms is significantly shorter than for the central and distributed scenarios. Indicating that in these scenarios short routes can be found when the demand is spatially clustered. This effect is however not present for the tight time window scenarios where modules are utilized (compare Figure \ref{fig:Modules_aggregated}). This shows how the utilization of service depots and multi-purpose vehicles changes the order in which requests are served and where vehicles are operating. 

\subsubsection{Impacts of module usage}

\begin{figure*}[ht]
	\centering
	\includegraphics[width=0.75\textwidth]{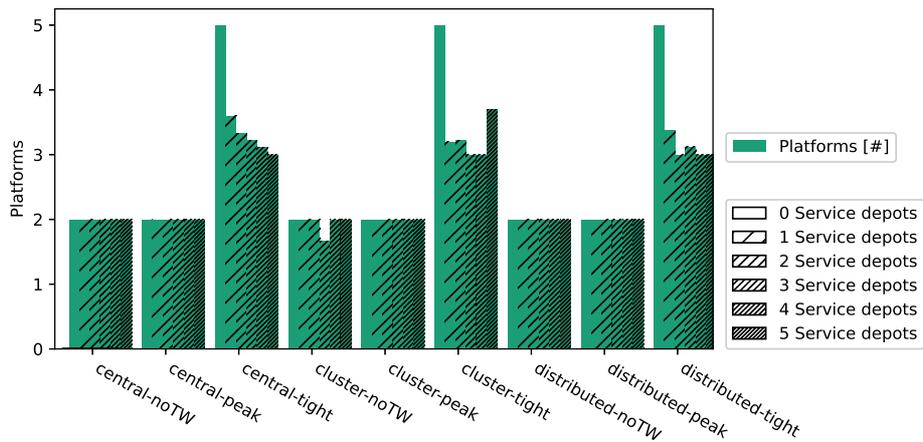}
	\caption{Average number of platforms in a solution for different number of service depots, spatial distributions and time window definitions.}
	\label{fig:Platforms_changes}
\end{figure*}

As described in the previous section, the number of platforms needed to serve the demand can be reduced by utilizing multi-purpose vehicles. In Figure \ref{fig:Platforms_changes} the number of platforms in a solution is shown with respect to the number of available service depots. Two main observations can be made from this figure. First, the visit of just one service depot, i.e. one module change, results in a reduction of approximately two platforms for all three spatial variations (central, cluster, and distributed). This is only true for scenarios with tight time window definitions. For other time window definitions, no or only small reduction in platforms due to visits at the service depot is observed. Second, no additional platforms are saved by having more service depots available. This indicates that there is an optimal combination of platforms and the number of modules, which cannot be improved by adding more service depots to the scenario.

\begin{figure*}[ht]
	\centering
	\includegraphics[width=0.85\textwidth]{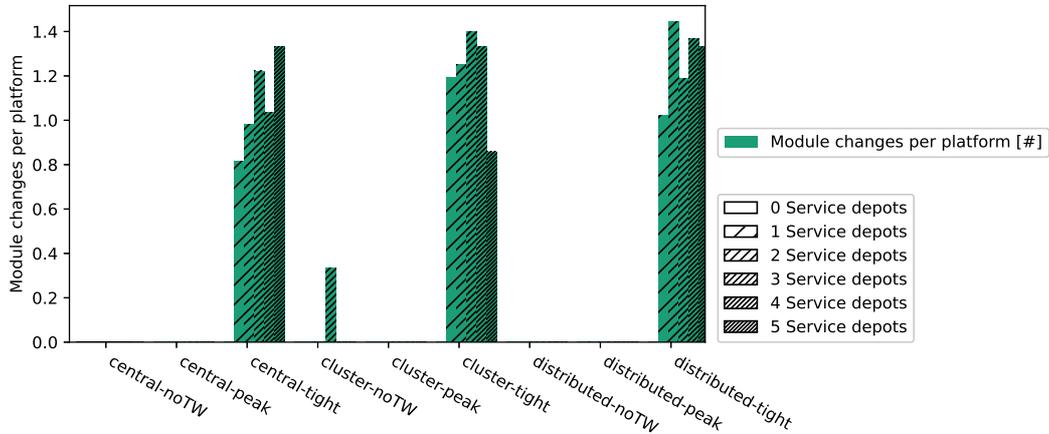}
	\caption{Module changes per platform for different number of service depots, averaged for different spatial distributions and time window definitions.}
	\label{fig:Module_changes_per_platform_scenario}
\end{figure*}

This observation is confirmed by the average number of module changes per platform shown in Figure \ref{fig:Module_changes_per_platform_scenario}. The different patterns in Figure \ref{fig:Module_changes_per_platform_scenario} indicate the different number of service depots in each scenario (compare Figure \ref{fig:scenarios}). Note that the number of service depots indicates the order in which they are added. Hence, the first service depot added is always service depot 1, the second is always service depot 2 and so on until all 5 service depots are added to the scenario. As described in the previous section, the results for scenarios with peak time windows and without time window constraints exhibit similar trends. Most of the platforms in these scenarios do not change their purpose along their route, meaning that not every platform in these scenarios visits a service depot along its route. This changes with tight time window definitions, where an average of 1.18 module changes are made per platform. This means that few platforms visit multiple service depots during their trip. It can also be seen that the number of module changes does not correlate with the number of available service depots, i.e. for an increasing number of available service depots, the number of module changes does not increase significantly. This suggests that several service depots are not or only rarely used and shows the importance of positioning the service depots at strategic locations within the service area.
\section{Discussion and Conclusion}\label{sec:conclusion}

In this study, a new vehicle technology for urban logistics and urban passenger transport is examined with a focus on the sequential consolidation of different passenger and freight flows in urban environments. By exchanging a removable module from the vehicle platform, different types of requests can be served using the same vehicle platform at different times. The paper proposes an extension of the pickup and delivery problem to model the novel operation of the said modular vehicles. We demonstrate the ability of the heuristic ALNS optimization algorithm to solve the routing problem for large-scale scenarios. 

The conducted experiments based on scenarios in Stockholm, Sweden showed that the use of multi-purpose vehicles lead to an objective value reduction of 11.3\%, 13.4\%, and 15.4\%, for scenarios with tight time windows and central, cluster and demand distribution, respectively. This reduction comes mainly from an approximate 35\% reduction in the number of platforms required to serve requests in these scenarios. This reduction comes at the cost of approximately 37\% more used modules. Furthermore, the results show that the use of service depots leads to a change in in which order the requests are served, resulting in an average reduction in the total duration of the vehicle trip of approximately 33\% and a reduction in travel distance of approximately 16\% for these scenarios. The number of unserved request remains unchanged when changing the vehicle technology, indicating an unchanged level of service.

Considering the usage of service depots, the average number of module changes is around 1.18 for scenarios with tight time window definitions, showing that every platform is in visiting at least one service depot along its route. Interestingly the number of platforms saved does not correlate with the number of service depot available, meaning that more available service depots and/or more service depot visits do not lead to a further reduction in required platforms. The main benefits are already achieved with a single visit. Lastly, it could be noted that the number of module changes is not significantly increased by increasing the number of available service depots in the scenario. This allows the conclusion that most service depots are not or only rarely used, and hence the positioning of service depots is crucial for the benefits of the proposed system. 

The results of this study can be used by practitioners and policymakers to conclude whether the combination of passenger and freight demand flows in each system will yield benefits compared to existing systems. In general, it can be said that for realistic (tight) time window definitions the new vehicle technology leads to promising results, mainly a reduction in required platforms to serve the same demand. Generally, in comparison with simultaneous integration concepts the proposed sequential integration concept guarantees that passenger requests are not affected by freight delivery processes. Hence, the level of service for passengers is higher in this form of integration.

Our study is subject to several limitations affecting the scalability of the model and the generality of the results reported. The computation time of the ALNS is capable of solving problems of around 200-400 requests within a reasonable time, however, for larger problems, the proposed heuristic algorithm needs to be improved with e.g. pre-clustering steps or multi-level optimization approaches. Further research may examine the sensitivity of numerical solutions to vehicle characteristics such as speed, range, and capacity. Additionally, the location and costs structures of depots and vehicle operations influence the generality of the proposed results. First, we assume that during every service depot visit the correct module type is always available. In the model this is guaranteed by setting the number of available modules per type and per service depot equal to the number of maximum allowed service depot visits, as defined by parameter $\vartheta$. In this way, we can guarantee successful operation of the transportation system, but we overestimate the module availability compared to a more realistic scenario. Second, the current version of the proposed model does not guarantee that the number of modules and their types is the same before and after the planning period at each depot. In practical operations this could require re-balancing of available vehicle fleets, which would reduce the benefits of the proposed multi-purpose transport system. Notwithstanding, the proposed modeling framework and optimization algorithm can be adapted to accommodate such considerations if necessary.

The proposed research can be extended in three directions. First, by integrating mixed-fleet properties (i.e., different vehicle sizes, simultaneous optimization of conventional and multi-purpose vehicles) for urban freight transport studies. Second, module scheduling and/or module inventory features can be added to the proposed three-index formulation. Third, the (service) depot positioning problem can be integrated in the proposed model in order to expand to strategic planning decisions.

\section{Acknowledgments}\label{sec:acknowledgements}
This work was in part supported by the Sustainable and Integrated Transport Systems (Hållbara och Integrerade Transport System HITS) project under grant number 2020-00565 (Vinnova). The computations were enabled by resources provided by the Swedish National Infrastructure for Computing (SNIC) at HPC2N partially funded by the Swedish Research Council through grant agreement no. 2018-05973.

\bibliography{references}

\begin{thebibliography}{53}
\expandafter\ifx\csname natexlab\endcsname\relax\def\natexlab#1{#1}\fi
\providecommand{\url}[1]{\texttt{#1}}
\providecommand{\href}[2]{#2}
\providecommand{\path}[1]{#1}
\providecommand{\DOIprefix}{doi:}
\providecommand{\ArXivprefix}{arXiv:}
\providecommand{\URLprefix}{URL: }
\providecommand{\Pubmedprefix}{pmid:}
\providecommand{\doi}[1]{\href{http://dx.doi.org/#1}{\path{#1}}}
\providecommand{\Pubmed}[1]{\href{pmid:#1}{\path{#1}}}
\providecommand{\bibinfo}[2]{#2}
\ifx\xfnm\relax \def\xfnm[#1]{\unskip,\space#1}\fi
\bibitem[{Alizadeh~Foroutan et~al.(2020)Alizadeh~Foroutan, Rezaeian and
  Mahdavi}]{alizadehforoutanGreenVehicleRouting2020}
\bibinfo{author}{Alizadeh~Foroutan, R.}, \bibinfo{author}{Rezaeian, J.},
  \bibinfo{author}{Mahdavi, I.}, \bibinfo{year}{2020}.
\newblock \bibinfo{title}{Green vehicle routing and scheduling problem with
  heterogeneous fleet including reverse logistics in the form of collecting
  returned goods}.
\newblock \bibinfo{journal}{Applied Soft Computing} \bibinfo{volume}{94},
  \bibinfo{pages}{106462}.
\newblock \DOIprefix\doi{10.1016/j.asoc.2020.106462}.
\bibitem[{Arslan et~al.(2016)Arslan, Agatz, Kroon and
  Zuidwijk}]{arslanCrowdsourcedDeliveryDynamic2016}
\bibinfo{author}{Arslan, A.}, \bibinfo{author}{Agatz, N.},
  \bibinfo{author}{Kroon, L.}, \bibinfo{author}{Zuidwijk, R.},
  \bibinfo{year}{2016}.
\newblock \bibinfo{title}{Crowdsourced {{Delivery}}: {{A Dynamic Pickup}} and
  {{Delivery Problem}} with {{Ad}}-{{Hoc Drivers}}}.
\newblock \bibinfo{type}{{{SSRN Scholarly Paper}}} \bibinfo{number}{ID
  2726731}. {Social Science Research Network}. \bibinfo{address}{{Rochester,
  NY}}.
\newblock \DOIprefix\doi{10.2139/ssrn.2726731}.
\bibitem[{Behiri et~al.(2018)Behiri, {Belmokhtar-Berraf} and
  Chu}]{behiriUrbanFreightTransport2018a}
\bibinfo{author}{Behiri, W.}, \bibinfo{author}{{Belmokhtar-Berraf}, S.},
  \bibinfo{author}{Chu, C.}, \bibinfo{year}{2018}.
\newblock \bibinfo{title}{Urban freight transport using passenger rail network:
  {{Scientific}} issues and quantitative analysis}.
\newblock \bibinfo{journal}{Transportation Research Part E: Logistics and
  Transportation Review} \bibinfo{volume}{115}, \bibinfo{pages}{227--245}.
\newblock \DOIprefix\doi{10.1016/j.tre.2018.05.002}.
\bibitem[{Beirigo et~al.(2018)Beirigo, Schulte and
  Negenborn}]{beirigoIntegratingPeopleFreight2018}
\bibinfo{author}{Beirigo, B.A.}, \bibinfo{author}{Schulte, F.},
  \bibinfo{author}{Negenborn, R.R.}, \bibinfo{year}{2018}.
\newblock \bibinfo{title}{Integrating {{People}} and {{Freight Transportation
  Using Shared Autonomous Vehicles}} with {{Compartments}}}.
\newblock \bibinfo{journal}{IFAC-PapersOnLine} \bibinfo{volume}{51},
  \bibinfo{pages}{392--397}.
\newblock \DOIprefix\doi{10.1016/j.ifacol.2018.07.064}.
\bibitem[{Berbeglia et~al.(2010)Berbeglia, Cordeau and
  Laporte}]{berbegliaDynamicPickupDelivery2010}
\bibinfo{author}{Berbeglia, G.}, \bibinfo{author}{Cordeau, J.F.},
  \bibinfo{author}{Laporte, G.}, \bibinfo{year}{2010}.
\newblock \bibinfo{title}{Dynamic pickup and delivery problems}.
\newblock \bibinfo{journal}{European Journal of Operational Research}
  \bibinfo{volume}{202}, \bibinfo{pages}{8--15}.
\newblock \DOIprefix\doi{10.1016/j.ejor.2009.04.024}.
\bibitem[{B{\"o}rjesson and
  Eliasson(2014)}]{borjessonExperiencesSwedishValue2014}
\bibinfo{author}{B{\"o}rjesson, M.}, \bibinfo{author}{Eliasson, J.},
  \bibinfo{year}{2014}.
\newblock \bibinfo{title}{Experiences from the {{Swedish Value}} of {{Time}}
  study}.
\newblock \bibinfo{journal}{Transportation Research Part A: Policy and
  Practice} \bibinfo{volume}{59}, \bibinfo{pages}{144--158}.
\newblock \DOIprefix\doi{10.1016/j.tra.2013.10.022}.
\bibitem[{Ceder and Wilson(1986)}]{cederBusNetworkDesign1986}
\bibinfo{author}{Ceder, A.}, \bibinfo{author}{Wilson, N.H.},
  \bibinfo{year}{1986}.
\newblock \bibinfo{title}{Bus network design}.
\newblock \bibinfo{journal}{Transportation Research Part B: Methodological}
  \bibinfo{volume}{20}, \bibinfo{pages}{331--344}.
\bibitem[{Chew et~al.(2013)Chew, Lee and
  Seow}]{chewGeneticAlgorithmBiobjective2013}
\bibinfo{author}{Chew, J.S.C.}, \bibinfo{author}{Lee, L.S.},
  \bibinfo{author}{Seow, H.V.}, \bibinfo{year}{2013}.
\newblock \bibinfo{title}{Genetic {{Algorithm}} for {{Biobjective Urban Transit
  Routing Problem}}}.
\newblock \bibinfo{journal}{Journal of Applied Mathematics}
  \bibinfo{volume}{2013}, \bibinfo{pages}{1--15}.
\newblock \DOIprefix\doi{10.1155/2013/698645}.
\bibitem[{Cochrane et~al.(2017)Cochrane, Saxe, Roorda and
  Shalaby}]{cochraneMovingFreightPublic2017}
\bibinfo{author}{Cochrane, K.}, \bibinfo{author}{Saxe, S.},
  \bibinfo{author}{Roorda, M.J.}, \bibinfo{author}{Shalaby, A.},
  \bibinfo{year}{2017}.
\newblock \bibinfo{title}{Moving freight on public transit: {{Best}} practices,
  challenges, and opportunities}.
\newblock \bibinfo{journal}{International Journal of Sustainable
  Transportation} \bibinfo{volume}{11}, \bibinfo{pages}{120--132}.
\newblock \DOIprefix\doi{10.1080/15568318.2016.1197349}.
\bibitem[{Cordeau and Laporte(2007)}]{cordeauDialarideProblemModels2007}
\bibinfo{author}{Cordeau, J.F.}, \bibinfo{author}{Laporte, G.},
  \bibinfo{year}{2007}.
\newblock \bibinfo{title}{The dial-a-ride problem: Models and algorithms}.
\newblock \bibinfo{journal}{Annals of Operations Research}
  \bibinfo{volume}{153}, \bibinfo{pages}{29--46}.
\newblock \DOIprefix\doi{10.1007/s10479-007-0170-8}.
\bibitem[{Dantzig and Ramser(1959)}]{dantzigTruckDispatchingProblem1959}
\bibinfo{author}{Dantzig, G.B.}, \bibinfo{author}{Ramser, J.H.},
  \bibinfo{year}{1959}.
\newblock \bibinfo{title}{The {{Truck Dispatching Problem}}}.
\newblock \bibinfo{journal}{Management Science} \bibinfo{volume}{6},
  \bibinfo{pages}{80--91}.
\newblock \DOIprefix\doi{10.1287/mnsc.6.1.80}.
\bibitem[{Derigs et~al.(2013)Derigs, Pullmann and
  Vogel}]{derigsTruckTrailerRouting2013}
\bibinfo{author}{Derigs, U.}, \bibinfo{author}{Pullmann, M.},
  \bibinfo{author}{Vogel, U.}, \bibinfo{year}{2013}.
\newblock \bibinfo{title}{Truck and trailer routing\textemdash{{Problems}},
  heuristics and computational experience}.
\newblock \bibinfo{journal}{Computers \& Operations Research}
  \bibinfo{volume}{40}, \bibinfo{pages}{536--546}.
\newblock \DOIprefix\doi{10.1016/j.cor.2012.08.007}.
\bibitem[{Desrochers et~al.(1988)Desrochers, Lenstra, Savelsbergh and
  Soumis}]{PDP1988}
\bibinfo{author}{Desrochers, M.}, \bibinfo{author}{Lenstra, J.},
  \bibinfo{author}{Savelsbergh, M.}, \bibinfo{author}{Soumis, F.},
  \bibinfo{year}{1988}.
\newblock \bibinfo{title}{Vehicle routing with time windows: Optimization and
  approximation}.
\newblock \bibinfo{journal}{Veh Rout Methods Stud} \bibinfo{volume}{16}.
\bibitem[{DLR(2021)}]{UShiftIIDemonstrator}
\bibinfo{author}{DLR}, \bibinfo{year}{2021}.
\newblock \bibinfo{title}{U-{{Shift II}} ({{Demonstrator}})}.
\newblock
  \bibinfo{howpublished}{https://verkehrsforschung.dlr.de/en/projects/u-shift/u-shift-ii-demonstrator}.
\bibitem[{Drexl(2013)}]{drexlApplicationsVehicleRouting2013}
\bibinfo{author}{Drexl, M.}, \bibinfo{year}{2013}.
\newblock \bibinfo{title}{Applications of the vehicle routing problem with
  trailers and transshipments}.
\newblock \bibinfo{journal}{European Journal of Operational Research}
  \bibinfo{volume}{227}, \bibinfo{pages}{275--283}.
\newblock \DOIprefix\doi{10.1016/j.ejor.2012.12.015}.
\bibitem[{D{\"u}ndar et~al.(2021)D{\"u}ndar, {\"O}m{\"u}rg{\"o}n{\"u}l{\c s}en
  and Soysal}]{dundarReviewSustainableUrban2021}
\bibinfo{author}{D{\"u}ndar, H.}, \bibinfo{author}{{\"O}m{\"u}rg{\"o}n{\"u}l{\c
  s}en, M.}, \bibinfo{author}{Soysal, M.}, \bibinfo{year}{2021}.
\newblock \bibinfo{title}{A review on sustainable urban vehicle routing}.
\newblock \bibinfo{journal}{Journal of Cleaner Production}
  \bibinfo{volume}{285}, \bibinfo{pages}{125444}.
\newblock \DOIprefix\doi{10.1016/j.jclepro.2020.125444}.
\bibitem[{Ghilas et~al.(2016a)Ghilas, Demir and
  Van~Woensel}]{ghilasAdaptiveLargeNeighborhood2016}
\bibinfo{author}{Ghilas, V.}, \bibinfo{author}{Demir, E.},
  \bibinfo{author}{Van~Woensel, T.}, \bibinfo{year}{2016}a.
\newblock \bibinfo{title}{An adaptive large neighborhood search heuristic for
  the {{Pickup}} and {{Delivery Problem}} with {{Time Windows}} and {{Scheduled
  Lines}}}.
\newblock \bibinfo{journal}{Computers \& Operations Research}
  \bibinfo{volume}{72}, \bibinfo{pages}{12--30}.
\newblock \DOIprefix\doi{10.1016/j.cor.2016.01.018}.
\bibitem[{Ghilas et~al.(2016b)Ghilas, Demir and
  Woensel}]{ghilasPickupDeliveryProblem2016}
\bibinfo{author}{Ghilas, V.}, \bibinfo{author}{Demir, E.},
  \bibinfo{author}{Woensel, T.V.}, \bibinfo{year}{2016}b.
\newblock \bibinfo{title}{The pickup and delivery problem with time windows and
  scheduled lines}.
\newblock \bibinfo{journal}{INFOR: Information Systems and Operational
  Research} \bibinfo{volume}{54}, \bibinfo{pages}{147--167}.
\newblock \DOIprefix\doi{10.1080/03155986.2016.1166793}.
\bibitem[{Huber and Geiger(2017)}]{huberOrderMattersVariable2017}
\bibinfo{author}{Huber, S.}, \bibinfo{author}{Geiger, M.J.},
  \bibinfo{year}{2017}.
\newblock \bibinfo{title}{Order matters \textendash{} {{A Variable Neighborhood
  Search}} for the {{Swap}}-{{Body Vehicle Routing Problem}}}.
\newblock \bibinfo{journal}{European Journal of Operational Research}
  \bibinfo{volume}{263}, \bibinfo{pages}{419--445}.
\newblock \DOIprefix\doi{10.1016/j.ejor.2017.04.046}.
\bibitem[{Ko{\c c} et~al.(2016)Ko{\c c}, Bekta{\c s}, Jabali and
  Laporte}]{kocThirtyYearsHeterogeneous2016}
\bibinfo{author}{Ko{\c c}, {\c C}.}, \bibinfo{author}{Bekta{\c s}, T.},
  \bibinfo{author}{Jabali, O.}, \bibinfo{author}{Laporte, G.},
  \bibinfo{year}{2016}.
\newblock \bibinfo{title}{Thirty years of heterogeneous vehicle routing}.
\newblock \bibinfo{journal}{European Journal of Operational Research}
  \bibinfo{volume}{249}, \bibinfo{pages}{1--21}.
\newblock \DOIprefix\doi{10.1016/j.ejor.2015.07.020}.
\bibitem[{Lenstra and Kan(1981)}]{lenstraComplexityVehicleRouting22}
\bibinfo{author}{Lenstra, J.K.}, \bibinfo{author}{Kan, A.H.G.R.},
  \bibinfo{year}{1981}.
\newblock \bibinfo{title}{Complexity of vehicle routing and scheduling
  problems}.
\newblock \bibinfo{journal}{Networks} \bibinfo{volume}{11},
  \bibinfo{pages}{221--227}.
\newblock \DOIprefix\doi{10.1002/net.3230110211}.
\bibitem[{Li et~al.(2014)Li, Krushinsky, Reijers and {van
  Woensel}}]{liShareaRideProblemPeople2014}
\bibinfo{author}{Li, B.}, \bibinfo{author}{Krushinsky, D.},
  \bibinfo{author}{Reijers, H.A.}, \bibinfo{author}{{van Woensel}, T.},
  \bibinfo{year}{2014}.
\newblock \bibinfo{title}{The {{Share}}-a-{{Ride Problem}}: {{People}} and
  parcels sharing taxis}.
\newblock \bibinfo{journal}{European Journal of Operational Research}
  \bibinfo{volume}{238}, \bibinfo{pages}{31--40}.
\newblock \DOIprefix\doi{10.1016/j.ejor.2014.03.003}.
\bibitem[{Li et~al.(2016)Li, Krushinsky, Van~Woensel and
  Reijers}]{liAdaptiveLargeNeighborhood2016}
\bibinfo{author}{Li, B.}, \bibinfo{author}{Krushinsky, D.},
  \bibinfo{author}{Van~Woensel, T.}, \bibinfo{author}{Reijers, H.A.},
  \bibinfo{year}{2016}.
\newblock \bibinfo{title}{An adaptive large neighborhood search heuristic for
  the share-a-ride problem}.
\newblock \bibinfo{journal}{Computers \& Operations Research}
  \bibinfo{volume}{66}, \bibinfo{pages}{170--180}.
\newblock \DOIprefix\doi{10.1016/j.cor.2015.08.008}.
\bibitem[{Li et~al.(2021)Li, Shalaby, Roorda and Mao}]{liUrbanRailService2021}
\bibinfo{author}{Li, Z.}, \bibinfo{author}{Shalaby, A.},
  \bibinfo{author}{Roorda, M.J.}, \bibinfo{author}{Mao, B.},
  \bibinfo{year}{2021}.
\newblock \bibinfo{title}{Urban rail service design for collaborative passenger
  and freight transport}.
\newblock \bibinfo{journal}{Transportation Research Part E: Logistics and
  Transportation Review} \bibinfo{volume}{147}, \bibinfo{pages}{102205}.
\newblock \DOIprefix\doi{10.1016/j.tre.2020.102205}.
\bibitem[{Los et~al.(2020)Los, Schulte, Spaan and
  Negenborn}]{losValueInformationSharing2020}
\bibinfo{author}{Los, J.}, \bibinfo{author}{Schulte, F.},
  \bibinfo{author}{Spaan, M.T.J.}, \bibinfo{author}{Negenborn, R.R.},
  \bibinfo{year}{2020}.
\newblock \bibinfo{title}{The value of information sharing for platform-based
  collaborative vehicle routing}.
\newblock \bibinfo{journal}{Transportation Research Part E: Logistics and
  Transportation Review} \bibinfo{volume}{141}, \bibinfo{pages}{102011}.
\newblock \DOIprefix\doi{10.1016/j.tre.2020.102011}.
\bibitem[{Marinov et~al.(2013)Marinov, Giubilei, Gerhardt, {\"O}zkan, Stergiou,
  Papadopol and Cabecinha}]{marinovUrbanFreightMovement2013}
\bibinfo{author}{Marinov, M.}, \bibinfo{author}{Giubilei, F.},
  \bibinfo{author}{Gerhardt, M.}, \bibinfo{author}{{\"O}zkan, T.},
  \bibinfo{author}{Stergiou, E.}, \bibinfo{author}{Papadopol, M.},
  \bibinfo{author}{Cabecinha, L.}, \bibinfo{year}{2013}.
\newblock \bibinfo{title}{Urban freight movement by rail}.
\newblock \bibinfo{journal}{Journal of Transport Literature}
  \bibinfo{volume}{7}, \bibinfo{pages}{87--116}.
\newblock \DOIprefix\doi{10.1590/S2238-10312013000300005}.
\bibitem[{Masson et~al.(2013)Masson, Lehu{\'e}d{\'e} and
  P{\'e}ton}]{massonAdaptiveLargeNeighborhood2013}
\bibinfo{author}{Masson, R.}, \bibinfo{author}{Lehu{\'e}d{\'e}, F.},
  \bibinfo{author}{P{\'e}ton, O.}, \bibinfo{year}{2013}.
\newblock \bibinfo{title}{An {{Adaptive Large Neighborhood Search}} for the
  {{Pickup}} and {{Delivery Problem}} with {{Transfers}}}.
\newblock \bibinfo{journal}{Transportation Science} \bibinfo{volume}{47},
  \bibinfo{pages}{344--355}.
\newblock \DOIprefix\doi{10.1287/trsc.1120.0432}.
\bibitem[{Melachrinoudis et~al.(2007)Melachrinoudis, Ilhan and
  Min}]{melachrinoudisDialarideProblemClient2007}
\bibinfo{author}{Melachrinoudis, E.}, \bibinfo{author}{Ilhan, A.B.},
  \bibinfo{author}{Min, H.}, \bibinfo{year}{2007}.
\newblock \bibinfo{title}{A dial-a-ride problem for client transportation in a
  health-care organization}.
\newblock \bibinfo{journal}{Computers \& Operations Research}
  \bibinfo{volume}{34}, \bibinfo{pages}{742--759}.
\newblock \DOIprefix\doi{10.1016/j.cor.2005.03.024}.
\bibitem[{Milit{\~a}o and Tirachini(2021)}]{militaoOptimalFleetSize2021}
\bibinfo{author}{Milit{\~a}o, A.M.}, \bibinfo{author}{Tirachini, A.},
  \bibinfo{year}{2021}.
\newblock \bibinfo{title}{Optimal fleet size for a shared demand-responsive
  transport system with human-driven vs automated vehicles: {{A}} total cost
  minimization approach}.
\newblock \bibinfo{journal}{Transportation Research Part A: Policy and
  Practice} \bibinfo{volume}{151}, \bibinfo{pages}{52--80}.
\newblock \DOIprefix\doi{10.1016/j.tra.2021.07.004}.
\bibitem[{Mourad et~al.(2019)Mourad, Puchinger and
  Chu}]{mouradSurveyModelsAlgorithms2019}
\bibinfo{author}{Mourad, A.}, \bibinfo{author}{Puchinger, J.},
  \bibinfo{author}{Chu, C.}, \bibinfo{year}{2019}.
\newblock \bibinfo{title}{A survey of models and algorithms for optimizing
  shared mobility}.
\newblock \bibinfo{journal}{Transportation Research Part B: Methodological}
  \bibinfo{volume}{123}, \bibinfo{pages}{323--346}.
\newblock \DOIprefix\doi{10.1016/j.trb.2019.02.003}.
\bibitem[{Mourad et~al.(2020)Mourad, Puchinger and
  Woensel}]{mouradIntegratingAutonomousDelivery2020}
\bibinfo{author}{Mourad, A.}, \bibinfo{author}{Puchinger, J.},
  \bibinfo{author}{Woensel, T.V.}, \bibinfo{year}{2020}.
\newblock \bibinfo{title}{Integrating autonomous delivery service into a
  passenger transportation system}.
\newblock \bibinfo{journal}{International Journal of Production Research}
  \bibinfo{volume}{0}, \bibinfo{pages}{1--24}.
\newblock \DOIprefix\doi{10.1080/00207543.2020.1746850}.
\bibitem[{{Mueller-Eberstein} and Franke(2000)}]{article}
\bibinfo{author}{{Mueller-Eberstein}, F.}, \bibinfo{author}{Franke, M.},
  \bibinfo{year}{2000}.
\newblock \bibinfo{title}{Project {{CarGo}} tram}.
\newblock
  \bibinfo{journal}{https://www.researchgate.net/publication/297113963\_Project\_CarGo\_Tram/}
  \bibinfo{volume}{124}, \bibinfo{pages}{337--340}.
\bibitem[{Ozturk and Patrick(2018)}]{ozturkOptimizationModelFreight2018}
\bibinfo{author}{Ozturk, O.}, \bibinfo{author}{Patrick, J.},
  \bibinfo{year}{2018}.
\newblock \bibinfo{title}{An optimization model for freight transport using
  urban rail transit}.
\newblock \bibinfo{journal}{European Journal of Operational Research}
  \bibinfo{volume}{267}, \bibinfo{pages}{1110--1121}.
\newblock \DOIprefix\doi{10.1016/j.ejor.2017.12.010}.
\bibitem[{Parragh(2011)}]{parraghIntroducingHeterogeneousUsers2011}
\bibinfo{author}{Parragh, S.N.}, \bibinfo{year}{2011}.
\newblock \bibinfo{title}{Introducing heterogeneous users and vehicles into
  models and algorithms for the dial-a-ride problem}.
\newblock \bibinfo{journal}{Transportation Research Part C: Emerging
  Technologies} \bibinfo{volume}{19}, \bibinfo{pages}{912--930}.
\newblock \DOIprefix\doi{10.1016/j.trc.2010.06.002}.
\bibitem[{Parragh and Cordeau(2017)}]{parraghBranchandpriceAdaptiveLarge2017}
\bibinfo{author}{Parragh, S.N.}, \bibinfo{author}{Cordeau, J.F.},
  \bibinfo{year}{2017}.
\newblock \bibinfo{title}{Branch-and-price and adaptive large neighborhood
  search for the truck and trailer routing problem with time windows}.
\newblock \bibinfo{journal}{Computers \& Operations Research}
  \bibinfo{volume}{83}, \bibinfo{pages}{28--44}.
\newblock \DOIprefix\doi{10.1016/j.cor.2017.01.020}.
\bibitem[{Parragh et~al.(2012)Parragh, Cordeau, Doerner and
  Hartl}]{parraghModelsAlgorithmsHeterogeneous2012}
\bibinfo{author}{Parragh, S.N.}, \bibinfo{author}{Cordeau, J.F.},
  \bibinfo{author}{Doerner, K.F.}, \bibinfo{author}{Hartl, R.F.},
  \bibinfo{year}{2012}.
\newblock \bibinfo{title}{Models and algorithms for the heterogeneous
  dial-a-ride problem with driver-related constraints}.
\newblock \bibinfo{journal}{OR Spectrum} \bibinfo{volume}{34},
  \bibinfo{pages}{593--633}.
\newblock \DOIprefix\doi{10.1007/s00291-010-0229-9}.
\bibitem[{Pisinger and Ropke(2007)}]{pisingerGeneralHeuristicVehicle2007}
\bibinfo{author}{Pisinger, D.}, \bibinfo{author}{Ropke, S.},
  \bibinfo{year}{2007}.
\newblock \bibinfo{title}{A general heuristic for vehicle routing problems}.
\newblock \bibinfo{journal}{Computers \& Operations Research}
  \bibinfo{volume}{34}, \bibinfo{pages}{2403--2435}.
\newblock \DOIprefix\doi{10.1016/j.cor.2005.09.012}.
\bibitem[{Pisinger and
  R{\o}pke(2010)}]{davidpisingerLargeNeighborhoodSearch2010}
\bibinfo{author}{Pisinger, D.}, \bibinfo{author}{R{\o}pke, S.},
  \bibinfo{year}{2010}.
\newblock \bibinfo{title}{Large {{Neighborhood Search}}}, in:
  \bibinfo{editor}{Gendreau, M.} (Ed.), \bibinfo{booktitle}{Handbook of
  {{Metaheuristics}}}. \bibinfo{edition}{second} ed..
  \bibinfo{publisher}{{Springer}}, p. \bibinfo{pages}{399420}.
\bibitem[{Qu and Bard()}]{qu_heterogeneous_2013}
\bibinfo{author}{Qu, Y.}, \bibinfo{author}{Bard, J.F.}, .
\newblock \bibinfo{title}{The heterogeneous pickup and delivery problem with
  configurable vehicle capacity}.
\newblock \bibinfo{journal}{Transportation Research Part C: Emerging
  Technologies} \bibinfo{volume}{32}, \bibinfo{pages}{1--20}.
\newblock \URLprefix
  \url{https://www.sciencedirect.com/science/article/pii/S0968090X13000612},
  \DOIprefix\doi{10.1016/j.trc.2013.03.007}.
\bibitem[{Rekiek et~al.(2006)Rekiek, Delchambre and
  Saleh}]{rekiekHandicappedPersonTransportation2006}
\bibinfo{author}{Rekiek, B.}, \bibinfo{author}{Delchambre, A.},
  \bibinfo{author}{Saleh, H.A.}, \bibinfo{year}{2006}.
\newblock \bibinfo{title}{Handicapped {{Person Transportation}}: {{An}}
  application of the {{Grouping Genetic Algorithm}}}.
\newblock \bibinfo{journal}{Engineering Applications of Artificial
  Intelligence} \bibinfo{volume}{19}, \bibinfo{pages}{511--520}.
\newblock \DOIprefix\doi{10.1016/j.engappai.2005.12.013}.
\bibitem[{Ropke and Pisinger(2006)}]{ropkeAdaptiveLargeNeighborhood2006}
\bibinfo{author}{Ropke, S.}, \bibinfo{author}{Pisinger, D.},
  \bibinfo{year}{2006}.
\newblock \bibinfo{title}{An {{Adaptive Large Neighborhood Search Heuristic}}
  for the {{Pickup}} and {{Delivery Problem}} with {{Time Windows}}}.
\newblock \bibinfo{journal}{Transportation Science} \bibinfo{volume}{40},
  \bibinfo{pages}{455--472}.
\newblock \DOIprefix\doi{10.1287/trsc.1050.0135}.
\bibitem[{Sacramento et~al.(2019)Sacramento, Pisinger and
  Ropke}]{sacramentoAdaptiveLargeNeighborhood2019}
\bibinfo{author}{Sacramento, D.}, \bibinfo{author}{Pisinger, D.},
  \bibinfo{author}{Ropke, S.}, \bibinfo{year}{2019}.
\newblock \bibinfo{title}{An adaptive large neighborhood search metaheuristic
  for the vehicle routing problem with drones}.
\newblock \bibinfo{journal}{Transportation Research Part C: Emerging
  Technologies} \bibinfo{volume}{102}, \bibinfo{pages}{289--315}.
\newblock \DOIprefix\doi{10.1016/j.trc.2019.02.018}.
\bibitem[{Savelsbergh and
  Woensel(2016)}]{savelsbergh50thAnniversaryInvited2016}
\bibinfo{author}{Savelsbergh, M.}, \bibinfo{author}{Woensel, T.V.},
  \bibinfo{year}{2016}.
\newblock \bibinfo{title}{50th {{Anniversary Invited Article}}\textemdash{{City
  Logistics}}: {{Challenges}} and {{Opportunities}}}.
\newblock \bibinfo{journal}{Transportation Science}
  \DOIprefix\doi{10.1287/trsc.2016.0675}.
\bibitem[{Savelsbergh and Sol(1995)}]{savelsberghGeneralPickupDelivery1995}
\bibinfo{author}{Savelsbergh, M.W.P.}, \bibinfo{author}{Sol, M.},
  \bibinfo{year}{1995}.
\newblock \bibinfo{title}{The {{General Pickup}} and {{Delivery Problem}}}.
\newblock \bibinfo{journal}{Transportation Science} \bibinfo{volume}{29},
  \bibinfo{pages}{17--29}.
\newblock \DOIprefix\doi{10.1287/trsc.29.1.17}.
\bibitem[{{Scania}(2020)}]{scaniaNXTConceptVehicle2020}
\bibinfo{author}{{Scania}}, \bibinfo{year}{2020}.
\newblock \bibinfo{title}{{{NXT}} concept vehicle represents a vision of the
  future}.
\newblock
  \bibinfo{howpublished}{https://www.scania.com/uk/en/home/experience-scania/features/nxt-concept-vehicle-represents-a-vision-of-the-future.html}.
\bibitem[{Schr{\"o}der and
  Liedtke(2017)}]{schroderIntegratedMultiagentUrban2017}
\bibinfo{author}{Schr{\"o}der, S.}, \bibinfo{author}{Liedtke, G.T.},
  \bibinfo{year}{2017}.
\newblock \bibinfo{title}{Towards an integrated multi-agent urban transport
  model of passenger and freight}.
\newblock \bibinfo{journal}{Research in Transportation Economics}
  \bibinfo{volume}{64}, \bibinfo{pages}{3--12}.
\newblock \DOIprefix\doi{10.1016/j.retrec.2016.12.001}.
\bibitem[{Shaw(1997)}]{shawNewLocalSearch1997a}
\bibinfo{author}{Shaw, P.}, \bibinfo{year}{1997}.
\newblock \bibinfo{title}{A {{New Local Search Algorithm Providing High Quality
  Solutions}} to {{Vehicle Routing Problems}}}.
\bibitem[{SteadieSeifi et~al.(2014)SteadieSeifi, Dellaert, Nuijten, Van~Woensel
  and Raoufi}]{steadieseifiMultimodalFreightTransportation2014}
\bibinfo{author}{SteadieSeifi, M.}, \bibinfo{author}{Dellaert, N.P.},
  \bibinfo{author}{Nuijten, W.}, \bibinfo{author}{Van~Woensel, T.},
  \bibinfo{author}{Raoufi, R.}, \bibinfo{year}{2014}.
\newblock \bibinfo{title}{Multimodal freight transportation planning: {{A}}
  literature review}.
\newblock \bibinfo{journal}{European Journal of Operational Research}
  \bibinfo{volume}{233}, \bibinfo{pages}{1--15}.
\newblock \DOIprefix\doi{10.1016/j.ejor.2013.06.055}.
\bibitem[{Tellez et~al.()Tellez, Vercraene, Lehuédé, Péton and
  Monteiro}]{tellez_fleet_2018}
\bibinfo{author}{Tellez, O.}, \bibinfo{author}{Vercraene, S.},
  \bibinfo{author}{Lehuédé, F.}, \bibinfo{author}{Péton, O.},
  \bibinfo{author}{Monteiro, T.}, .
\newblock \bibinfo{title}{The fleet size and mix dial-a-ride problem with
  reconfigurable vehicle capacity}.
\newblock \bibinfo{journal}{Transportation Research Part C: Emerging
  Technologies} \bibinfo{volume}{91}, \bibinfo{pages}{99--123}.
\newblock \URLprefix
  \url{https://www.sciencedirect.com/science/article/pii/S0968090X18303851},
  \DOIprefix\doi{10.1016/j.trc.2018.03.020}.
\bibitem[{Todosijevi{\'c} et~al.(2017)Todosijevi{\'c}, Hanafi, Uro{\v
  s}evi{\'c}, Jarboui and Gendron}]{todosijevicGeneralVariableNeighborhood2017}
\bibinfo{author}{Todosijevi{\'c}, R.}, \bibinfo{author}{Hanafi, S.},
  \bibinfo{author}{Uro{\v s}evi{\'c}, D.}, \bibinfo{author}{Jarboui, B.},
  \bibinfo{author}{Gendron, B.}, \bibinfo{year}{2017}.
\newblock \bibinfo{title}{A general variable neighborhood search for the
  swap-body vehicle routing problem}.
\newblock \bibinfo{journal}{Computers \& Operations Research}
  \bibinfo{volume}{78}, \bibinfo{pages}{468--479}.
\newblock \DOIprefix\doi{10.1016/j.cor.2016.01.016}.
\bibitem[{Toffolo et~al.(2018)Toffolo, Christiaens, {van Malderen}, Wauters and
  Vanden~Berghe}]{toffoloStochasticLocalSearch2018}
\bibinfo{author}{Toffolo, T.A.}, \bibinfo{author}{Christiaens, J.},
  \bibinfo{author}{{van Malderen}, S.}, \bibinfo{author}{Wauters, T.},
  \bibinfo{author}{Vanden~Berghe, G.}, \bibinfo{year}{2018}.
\newblock \bibinfo{title}{Stochastic local search with learning automaton for
  the swap-body vehicle routing problem}.
\newblock \bibinfo{journal}{Computers \& Operations Research}
  \bibinfo{volume}{89}, \bibinfo{pages}{68--81}.
\newblock \DOIprefix\doi{10.1016/j.cor.2017.08.002}.
\bibitem[{Toth and Vigo(2014)}]{tothVehicleRoutingProblems2014}
\bibinfo{editor}{Toth, P.}, \bibinfo{editor}{Vigo, D.} (Eds.),
  \bibinfo{year}{2014}.
\newblock \bibinfo{title}{Vehicle Routing: {{Problems}}, Methods, and
  Applications}.
\newblock {{MOS}}-{{SIAM}} Series on Optimization. \bibinfo{edition}{second
  edition} ed., \bibinfo{publisher}{{Society for Industrial and Applied
  Mathematics; Mathematical Optimization Society}},
  \bibinfo{address}{{Philadelphia}}.
\bibitem[{{VeRoLog}(2014)}]{verologVeRoLogSolverChallenge2014}
\bibinfo{author}{{VeRoLog}}, \bibinfo{year}{2014}.
\newblock \bibinfo{title}{{{VeRoLog Solver Challenge}} 2014 \textendash{}
  {{Verolog}}}.
\newblock
  \bibinfo{howpublished}{https://www.euro-online.org/websites/verolog/verolog-solver-challenge-2014/}.

\end{thebibliography}

\end{document}